\newcommand{\End}{\mathop{\mathrm{End}}\nolimits}
\newcommand{\Hom}{\mathop{\mathrm{Hom}}\nolimits}
\newcommand{\indlim}{\mathop{\underrightarrow{\lim}}\limits}
\newcommand{\overk}{\stackrel{/k}{\hookrightarrow}}
\newcommand{\Pic}{\mathop{\mathrm{Pic}}\nolimits}
\newcommand{\prlim}{\mathop{\underleftarrow{\lim}}\limits}
\newcommand{\Sm}{\mathop{\mathcal{S}\mathrm{m}}\nolimits}
\newcommand{\SOm}{\mathop{\underline{\Omega}}}
\documentclass[11pt]{amsart}
\usepackage[OT2,T2A]{fontenc}
\usepackage{amssymb} 
\advance \topmargin by -\headheight
\advance \topmargin by -\headsep
\newtheorem{theorem}{Theorem}[section]
\newtheorem{lemma}[theorem]{Lemma}

\newtheorem{proposition}[theorem]{Proposition}
\newtheorem{conjecture}[theorem]{Conjecture}
\newtheorem{conj-corollary}[theorem]{``Corollary''}
\newtheorem{cor}[theorem]{Corollary}
\title{Stable birational invariants with Galois descent and differential forms} 
\author{M.Rovinsky} 
\address{National Research University Higher School of Economics, 
Laboratory of Algebraic Geometry, 7 Vavilova Str., 
Moscow 117312, Russia \& Institute for Information Transmission Problems 
of Russian Academy of Sciences \& Independent University of Moscow}
\email{marat@mccme.ru}
\thanks{The author was supported by the National Science Foundation 
under agreement No. DMS-0635607. During the final write-up, the author 
was partially supported by RFBR grant 10-01-93113-CNRSL-a and by
AG Laboratory GU-HSE, RF government grant, ag. 11 11.G34.31.0023}
\date{}
\begin{document}
\begin{abstract} I show that the cohomology of the generic points of algebraic 
complex varieties becomes {\sl stable} birational invariant, when considered 
`modulo the cohomology of the generic points of the affine spaces'. 
\end{abstract} 
\maketitle

These notes are concerned with certain birational invariants of smooth 
algebraic varieties. All such invariants are {\sl dominant} sheaves, 
cf. below; the dominant sheaves are characterized in Proposition 
\ref{bir-A1-presheaf}. 

Two classes of invariants are of special interest: (i) {\sl stable}, 
i.e., taking the same values on a variety and on its direct product with 
an affine space, and (ii) constant on the projective spaces. Though the 
latter class is a priori wider, there are no known examples of non-stable 
invariants vanishing on the projective spaces. Here an attempt of 
comparison is made. Namely, it is shown that the corresponding adjoint 
functors coincide on the following types of invariants: (i) of `level 1', 
cf. Proposition \ref{pic-as-adj} and also p.\pageref{coniveau-filtr}, 
(ii) `related to cohomology' (or to closed differential forms). 

Differential forms play a very special r\^ole in the story, cf. e.g. 
Conjecture \ref{irred-diff-for}. Moreover, all known examples of 
{\sl simple} invariants (as objects of an abelian category) `come from' 
differential forms: except for two invariants related to the multiplicative 
and the additive groups ($Y\mapsto(k(Y)^{\times}/k^{\times})_{{\mathbb Q}}$ 
and $Y\mapsto k(Y)/k$, the logarithmic and the exact differentials, cf. 
below), they are values of the functor ${\mathbb B}^0$ from 
\S\ref{Irred-obj}. For these reasons the differential forms are studied in 
detail. It is shown in Corollary \ref{log-diff-form} that the cohomology of 
the generic points of algebraic (complex) varieties becomes {\sl stable} 
birational invariant, when considered `modulo the cohomology of the generic 
points of the affine spaces'. 

The principal new results of \S\ref{Differential-forms} are 
Propositions \ref{A1-quotient-of-closed} and \ref{1-forms}. It is shown in 
Proposition \ref{A1-quotient-of-closed} that (i) the quotient $V^{\bullet}$ 
of the sheaf of algebras of closed differential forms by the ideal generated 
by the exact 1-forms and the logarithmic differentials is stable and 
(ii) $V^{\bullet}$ is the maximal stable quotient of the sheaf of closed 
differential forms. Proposition \ref{1-forms} gives a complete description 
of the sheaf of closed 1-forms. 

Depending on what is more convenient, we shall consider our `invariants' 
either as dominant sheaves, or as representations, cf. \S\ref{Summ-equiv}. 
E.g., the simplicity is more natural in the context of representations. 

\section{Dominant presheaves and sheaves} \label{pre-sheaves} 
{\bf Notations.} From now on we fix an algebraically closed field 
$k$ of characteristic zero, and denote by $E$ a variable 
coefficient field of characteristic zero. 
Denote by $\text{Vec}_E$ the category of $E$-vector spaces. 

I am interested in birational invariants of (or ``presheaves on'') 
$k$-varieties. More precisely, let ${\mathcal S}m_k'$ be the category, 
whose objects are smooth $k$-varieties and the morphisms are smooth 
$k$-morphisms. Define the pretopology on ${\mathcal S}m_k'$ by saying 
that the covers are dominant morphisms. Recall, that a presheaf is a 
sheaf if the following diagram is an equalizer for any covering $Y\to X$: 
\begin{equation}\label{opr-puchka}{\mathcal F}(X)\to{\mathcal F}(Y)
\rightrightarrows{\mathcal F}(Y\times_XY).\end{equation} 

The category of the sheaves of $E$-vector spaces on this site 
is denoted by $\Sm_G(E)$ and $\Sm_G:=\Sm_G({\mathbb Q})$. 

{\it Example.} For each irreducible smooth $k$-variety $X$ and integer 
$0\le q\le\dim X$ let $\Psi_{X,q}:Y\mapsto Z^q(k(X)\otimes_kk(Y))
_{{\mathbb Q}}$ (${\mathbb Q}$-linear combinations of irreducible 
subvarieties on $X\times_kY$ of codimension $q$ dominant over $X$ and 
$Y$.) This is a sheaf. Set $\Psi_X:=\Psi_{X,\dim X}$. The sheaves 
$\Psi_X$ for all $X$ form a system of generators of $\Sm_G$. 

{\it Definition.} 1. A presheaf ${\mathcal F}$ is 
${\mathbb A}^1$-{\sl invariant} (or {\sl stable}) if 
${\mathcal F}(X)\stackrel{\sim}{\longrightarrow}
{\mathcal F}(X\times{\mathbb A}^1)$ for all $X$. 

2. Let ${\mathcal S}$ be a collection of dominant morphisms in 
${\mathcal S}m_k'$ with connected fibres. Assume that ${\mathcal S}$ is 
stable under base changes of its arbitrary element by itself: ${\rm pr}_1:
X\times_YX\to X$ belongs to ${\mathcal S}$ if $X\to Y$ belongs to 
${\mathcal S}$. A presheaf ${\mathcal F}$ is called an ${\mathcal S}$-presheaf 
if ${\mathcal F}(Y)\stackrel{\sim}{\longrightarrow}{\mathcal F}(X)$ 
for all $(X\to Y)\in{\mathcal S}$.  

Denote by $\Sm_G^{{\mathcal S}}$ the full subcategory in $\Sm_G$ 
consisting of ${\mathcal S}$-sheaves. More particularly, denote by 
${\mathcal I}_G$ the full subcategory in $\Sm_G$ consisting of 
${\mathbb A}^1$-invariant sheaves. Under assumptions of 
\S\ref{Properties-of}, $\Sm_G^{{\mathcal S}}\subseteq{\mathcal I}_G$. 

For any dominant presheaf ${\mathcal F}$ denote 
by $\underline{{\mathcal F}}$ its dominant sheafification. 

For each smooth $k$-variety $Y$, we denote by $\overline{Y}$ a smooth 
compactification of $Y$. 

\subsection{Examples of ${\mathbb A}^1$-invariant presheaves}
In this section we consider some examples of dominant presheaves with values 
in various abelian categories. They come either from algebro-geometric 
constructions, or from a cohomology theory $H^*$ (with 
coefficients in a commutative ${\mathbb Q}$-algebra $B$). As Example 5 
suggests, those of these examples that are ${\mathbb A}^1$-invariant sheaves, 
are related. This is one of motivations for Conjecture \ref{equiv-B}. 

An effective pure motive is a pair consisting of a smooth projective 
variety and a projector in the algebra of correspondences modulo 
numerical equivalence. Morphisms of {\sl co(ntra)variant} pure motives 
are defined by correspondences modulo numerical equivalence so that 
they behave as action on the (co)homology. 

Denote by ${\mathcal M}_k$ the category of covariant pure $k$-motives 
(and by ${\mathcal M}_k^{\text{op}}$ its opposite, the category of 
contravariant pure $k$-motives). By a well-known result of U.Jannsen, these 
two categories are abelian and semisimple. A simple effective pure motive 
is called {\sl primitive} if it is ``not divisible by the Lefschetz motive'', 
the motive $({\mathbb P}^1,\pi)$, where $\pi$ induces 0 on the 0-th and the 
identity on the second (co)homology. 

Denote by $\overline{Y}^{\text{prim}}$ the sum in the motive of 
$\overline{Y}$ of all its primitive submotives; $CH^q$ is the (Chow) 
group of codimension $q$ cycles modulo rational equivalence. 
We also use notations and identifications of \S\ref{obozn-predst}. 

\begin{tabular}{rlcc} &{\bf Invariant} of a connected $Y$ 
(dominant presheaf) & {\bf Values} &{\bf stable} \\
1&$K_q(Y)_{{\mathbb Q}}$ for $q\ge 0$/ its sheafification& 
$\text{Vec}_{{\mathbb Q}}$& yes/only for $q=0$ \\ 
2&$H^q(Y)$ for $q\ge 0$/ its sheafification $\underline{H}^q$& 
$B$-mod & yes/only for $q=0$ \\ 
3& $\Gamma(\overline{Y},\bigotimes^{\bullet}_{{\mathcal O}_{\overline{Y}}}
\Omega^1_{\overline{Y}|k})$ / its 
sheafification $\Gamma(\bigotimes^{\bullet}_F\Omega^1_{F|k})$, 
cf. Remark on p.\pageref{glok-por} & $\text{Vec}_k$ & yes/no \\ 
4&$\Phi^pCH^q(X\times_kk(Y))_{{\mathbb Q}}$ for a smooth $X$ and 
$\lefteqn{\text{a ``universal'' filtration $\Phi^{\bullet}$ on the 
Chow groups}}$ \\ 
&(e.g., $A(k(Y))_{{\mathbb Q}}$ for an abelian $k$-variety $A$)& 
$\text{Vec}_{{\mathbb Q}}$& yes \\ 
5&$\overline{Y}^{\text{prim}}=\bigoplus_M\overline{Y}^{\text{prim}}_M$ 
(multiplicity-one sheaf, by Proposition \ref{funktor-B})& 
${\mathcal M}_k^{\text{op}}$ & yes \\
6&$Z^{\dim Y}(F\otimes_kk(Y))_{{\mathbb Q}}$& 
$\Sm_G^{\text{op}}$& no \\
7&$Z^q(Y\times_kF)_{{\mathbb Q}}$ for $q\ge 0$ / its sheafification 
$Z^q(F\otimes_kk(Y))_{{\mathbb Q}}$& $\Sm_G$& only for $q=0$ \\ 
&(a) composition with the evaluation functor on $X$, i.e., $\Psi_{X,q}$& 
$\text{Vec}_{{\mathbb Q}}$ & only for $q=0$ \\
8&$C_{k(Y)}:={\mathcal I}\Psi_X$, cf. \S\ref{Properties-of}, 
(and its quotient $CH_0(\overline{Y}_F)_{{\mathbb Q}}$) 
& ${\mathcal I}_G^{\text{op}}$& yes \\ 
&(a) composition with $\Hom_{\Sm_G^{\text{op}}}(\underline{H}^q_c,-)$: \\
& $H^q(\overline{Y})/N^1=:\underline{H}^q_c(Y)$ for any $q\ge 0$ 
(subsheaf of the sheaf $\underline{H}^q$) & $B$-mod & yes \\ 
&(a$'$) $\lefteqn{\text{$k={\mathbb C}$: the image in 
$\underline{H}^{2q}(-({\mathbb C});{\mathbb Q})(Y)$ of the maximal Hodge 
substructure of $H^{2q}(\overline{Y}({\mathbb C}),{\mathbb Q})$ in}}$ \\ 
& $F^1$, cf. \S\ref{Differential-forms}, p.\pageref{Differential-forms}, 
(its vanishing is equivalent to the Hodge's conjecture) & 
$\text{Vec}_{{\mathbb Q}}$ & yes \\ 
&(b) composition with $\Hom_{\Sm_G^{\text{op}}}(\Omega_{F|k}^{\bullet},-)$: \\
& $\Gamma(\overline{Y},\Omega^{\bullet}_{\overline{Y}|k})$ (subsheaf of 
the sheaf $\underline{H}^{\bullet}_{\text{dR}|k,c}$)& $\text{Vec}_k$ & yes 
\end{tabular} 

Except for $\Gamma(\overline{Y},\bigotimes^{\bullet}
_{{\mathcal O}_{\overline{Y}}}\Omega^1_{\overline{Y}|k})$, 
all these invariants have Galois descent property. 
Except for $Z^q(Y\times_kF)_{{\mathbb Q}}$ for $q>0$, 
$K_q(Y)_{{\mathbb Q}}$ for $q\ge 0$ and $H^q(Y)$ for $q>0$, 
all these invariants are birational. ($N^1$ in example 8 (a) denotes 
the first term of the coniveau filtration on $H^*$.) 

Some of the above presheaves are defined using a compactification 
$\overline{Y}$. To show that each of such presheaves is in fact well-defined 
(and therefore, birationally invariant), one can use the facts that (i) 
any birational map is a composition of blow-ups and blow-downs with smooth 
centres, cf. \cite{akmw}, and (ii) the cohomology (resp., motive) of a 
blow-up is the direct sum of the cohomology of the original variety and 
of the Gysin image (resp., Tate twist) of the cohomology (resp., motive) 
of the subvariety which is blown up. Such a presheaf is 
${\mathbb A}^1$-invariant, since the cohomology (resp., motive) of the 
product of a proper variety $X$ with the projective line is the direct sum 
of the pull-back of the cohomology (resp., motive) of $X$ and of the Gysin 
image (resp., Tate twist) of the cohomology (resp., motive) 
of $X\times\{0\}\cong X$. 

To conclude that a birational ${\mathbb A}^1$-invariant presheaf 
is a sheaf, one checks that it has the Galois descent property, 
so Proposition \ref{bir-A1-presheaf} can be applied. 

\begin{lemma} \label{example-4} For an arbitrary commutative 
$k$-group $A$, let ${\mathcal H}^A_1$ be the presheaf 
$Y\mapsto\bigoplus_{y\in Y^0}(A(k(y))/A(k))_{{\mathbb Q}}$. Then 
${\mathcal H}^A_1$ is a sheaf; it is simple (=irreducible) for simple $A$. 
Let a presheaf ${\mathcal F}$ be the composition of the Picard 
functor $Y\mapsto\Pic^0(\overline{Y})$ with an additive functor 
on the category of abelian $k$-varieties, e.g. 
$\Pic^{\circ}_{{\mathbb Q}}:Y\mapsto\Pic^0(\overline{Y})_{{\mathbb Q}}$, 
$\underline{H}^1_c:Y\mapsto H^1(\overline{Y})$, 
or $\Omega^1_{|k,\text{{\rm reg}}}:Y\mapsto
\Gamma(\overline{Y},\Omega^1_{\overline{Y}|k})$. Then ${\mathcal F}$ 
is a sheaf and ${\mathcal F}=\bigoplus_A{\mathcal F}(\tilde{A})
\otimes_{\End(\tilde{A})}{\mathcal H}^{\tilde{A}}_1$, where $A$ runs through 
the isogeny classes of simple abelian $k$-varieties and $\tilde{A}$ 
is a representative of $A$. \end{lemma} 
Thus, such sheaves ${\mathcal F}$ are direct sums of copies 
of simple sheaves ${\mathcal H}^A_1$. 

{\it Proof.} Suppose that ${\mathcal A}$ is an abelian category. Then any 
semisimple object $N\in{\mathcal A}$ splits {\sl canonically} into the direct 
sum over the isomorphism classes $M$ of simple objects in ${\mathcal A}$ of 
its $M$-isotypical parts $N_M$. 
Clearly, for any representative $\tilde{M}$ of the isomorphism class $M$ the 
natural morphism 
$\Hom_{{\mathcal A}}(\tilde{M},N)\otimes_{\End(\tilde{M})}\tilde{M}\to 
N_M$ by $\varphi\otimes a\mapsto\varphi(a)$. This is an isomorphism. 
Applying an additive functor ${\mathfrak F}:{\mathcal A}\to
{\mathcal B}^{\text{op}}$ to the above isotypical decomposition of 
$N$, we get a canonical isomorphism $\prod_M{\mathfrak F}(\tilde{M})
\otimes_{\End(\tilde{M})}\Hom_{{\mathcal A}}(N,\tilde{M})
\stackrel{\sim}{\longrightarrow}{\mathfrak F}(N)$, 
$f\otimes l\mapsto l^{\ast}f$, where the following duality is used: 
$\Hom_{\text{mod-$\End(\tilde{M})$}}(\Hom_{{\mathcal A}}
(\tilde{M},N),\End(\tilde{M}))_{{\mathbb Q}}\stackrel{\sim}{\longrightarrow}
\Hom_{{\mathcal A}}(N,\tilde{M})$. (It is induced by the composition 
pairing $\Hom_{{\mathcal A}}(\tilde{M},N)\otimes
\Hom_{{\mathcal A}}(N,\tilde{M})\to\End(\tilde{M})$.)

As ${\mathcal A}$, we take either the category of abelian $k$-varieties 
with morphisms $\otimes{\mathbb Q}$, or the bigger category 
${\mathcal M}_k^{\text{op}}$. 
In the case of abelian varieties, the isomorphism classes of simple 
objects are the isogeny classes of simple abelian $k$-varieties, 
whereas the existence of the isotypical decomposition corresponds to 
the fact that for any abelian $k$-variety $B$ the natural morphism 
$\bigoplus_A\Hom_{\text{ab.$k$-var}}(\tilde{A},B)\otimes_{\End(\tilde{A})}
\tilde{A}\to B$, $\varphi\otimes a\mapsto\varphi(a)$, 
where $A$ runs through the isogeny classes of simple abelian $k$-varieties 
and $\tilde{A}$ is a representative of $A$, is an isogeny. 

Thus, any sheaf ${\mathcal F}$ with semisimple values in ${\mathcal A}$ also 
splits {\sl canonically} into the direct sum over the isomorphism classes $M$ 
of simple objects in ${\mathcal A}$ of its $M$-isotypical parts 
${\mathcal F}_M$. 

When ${\mathcal A}={\mathcal M}_k^{\text{op}}$ and ${\mathcal F}$ is 
the dominant sheaf $Y\mapsto\overline{Y}^{\text{prim}}$, we get that 
${\mathcal F}$ splits canonically into the direct sum of its $M$-isotypical 
parts $Y\mapsto\overline{Y}^{\text{prim}}_M$. By Proposition \ref{funktor-B}, 
the $M$-isotypical part $Y\mapsto\overline{Y}^{\text{prim}}_M$ is a simple 
sheaf. 

If $B={\rm Alb}(\overline{Y})$ (the Albanese variety) then 
$\Hom_{\text{ab.$k$-var}}(B,\tilde{A})=\tilde{A}(k(Y))/\tilde{A}(k)$, and thus, 
${\mathcal F}(Y)={\mathcal F}(B)\stackrel{\sim}{\longrightarrow}\bigoplus_A
{\mathcal F}(\tilde{A})\otimes_{\End(\tilde{A})}(\tilde{A}(k(Y))/\tilde{A}(k))$. 
It is quite evident that ${\mathcal H}^{\tilde{A}}_1$ is a sheaf. By 
Proposition \ref{bir-A1-presheaf}, in the case of abelian variety $\tilde{A}$, 
it suffices to check the Galois descent property, which is equivalent to the 
following one: for any abelian $k$-variety $\tilde{A}$ and any finite group 
$H$ of its automorphisms such that $H_0(H,\tilde{A})=0$ one has 
${\mathcal F}(\tilde{A})^H=0$. Clearly, this property holds. The simplicity 
of the sheaf ${\mathcal H}^{\tilde{A}}_1$ follows from the fact that for any 
algebraically closed field extension $K|k(\tilde{A})$ and for any subvariety 
$Z$ of $A$ of positive dimension there are no proper subgroups of 
$\tilde{A}(K)$ containing all generic $K$-points of $Z$. (Any point of 
$\tilde{A}$ is a sum of generic points of $\tilde{A}$; any sum of $\dim A$ 
generic $K$-points of $Z$ in sufficiently general position is a generic 
point of $Z$). This argument works more naturally in the context of 
representations, cf. \S\ref{Alternative-descr}. \qed 

\vspace{4mm}

{\it Remark.} For an abelian $k$-variety $A$, the sheaf 
$A_{{\mathbb Q}}:Y\mapsto A(k(Y))_{{\mathbb Q}}$ factors through the 
Albanese functor, but considered as a functor to the category of torsors 
over abelian $k$-varieties, so additive functors do not make sense 
and Lemma \ref{example-4} is not applicable to this sheaf. 
In particular, it is not semisimple. 

\vspace{4mm}

Propositions \ref{1-forms} and \ref{semisimpl-top} suggest that 
(i) the isomorphism classes of irreducible subquotients of 
$\underline{H}^{\bullet}_c$ are the same as that of 
$\Omega^{\bullet}_{|k,{\rm reg}}:Y\mapsto
\Gamma(\overline{Y},\Omega^{\bullet}_{\overline{Y}|k})$, 
(ii) they can be naturally identified with the irreducible 
effective primitive motives, and 
(iii) the isomorphism classes of irreducible subquotients 
of $\underline{H}^{\bullet}$ are related to more general irreducible 
effective motives, such as the Tate motive ${\mathbb Q}(-1)$ in the case 
of $\underline{H}^1_{{\rm dR}/k}$. 

\begin{lemma} Any dominant sheaf ${\mathcal F}$ with values in an abelian 
category with objects of finite length (e.g., in a category of 
finite-dimensional vector spaces) is ${\mathbb A}^1$-invariant. \end{lemma} 
{\it Proof.} Any smooth morphism of connected smooth $k$-varieties is covering, 
so $X\times({\mathbb A}^1\times{\mathbb A}^1\smallsetminus\Delta)\stackrel{p}
{\longrightarrow}X\times({\mathbb A}^1\times{\mathbb A}^1\smallsetminus\Delta)/
{\mathfrak S}_2$ is a cover for any $X$. On the other hand, it is the 
coequalizer of $X\times({\mathbb A}^1\times{\mathbb A}^1\smallsetminus\Delta)
\stackrel{id,(12)}\rightrightarrows 
X\times({\mathbb A}^1\times{\mathbb A}^1\smallsetminus\Delta)$. 
Therefore, ${\mathcal F}(X\times({\mathbb A}^1\times{\mathbb A}^1
\smallsetminus\Delta)/{\mathfrak S}_2)\stackrel{p^*}{\longrightarrow}
{\mathcal F}(X\times({\mathbb A}^1\times{\mathbb A}^1\smallsetminus\Delta))$ 
(i) is injective, (ii) factors through the ${\mathfrak S}_2$-invariants. 
As $({\mathbb A}^1\times{\mathbb A}^1\smallsetminus\Delta)/{\mathfrak S}_2
\cong{\mathbb A}^1\times{\mathbb A}^1\smallsetminus\Delta$($\cong
{\mathbb A}^1\times{\mathbb G}_m$), the source and the target of $p^*$ are 
isomorphic. As they are of finite length, the inclusion $p^*$ is an 
isomorphism. This implies that the involution $(12)$ is identical on 
${\mathcal F}(X\times({\mathbb A}^1\times{\mathbb A}^1\smallsetminus\Delta))$, 
so in the exact sequence, defining the sheaf condition for the cover 
$X\times{\mathbb A}^1\longrightarrow X$, $0\to{\mathcal F}(X)\to
{\mathcal F}(X\times{\mathbb A}^1)\rightrightarrows{\mathcal F}
(X\times{\mathbb A}^1\times{\mathbb A}^1)$ the double arrow consists 
of equal morphisms, i.e. ${\mathcal F}(X)\stackrel{\sim}{\longrightarrow}
{\mathcal F}(X\times{\mathbb A}^1)$. \qed 

\subsection{Properties of $\Sm_G^{{\mathcal S}}$} \label{Properties-of} 
Clearly, a subsheaf of an ${\mathcal S}$-sheaf is an ${\mathcal S}$-sheaf: 
if ${\mathcal G}$ is a subsheaf of an ${\mathcal S}$-sheaf ${\mathcal F}$ 
then for any $(Y\to X)\in{\mathcal S}$ the parallel arrows in the upper 
line in the commutative diagram
$$\begin{array}{ccccc}{\mathcal F}(X)&\longrightarrow&
{\mathcal F}(Y)&\rightrightarrows&{\mathcal F}(Y\times_XY)\\
\bigcup&&\bigcup&&\bigcup\\
{\mathcal G}(X)&\longrightarrow&{\mathcal G}(Y)&
\rightrightarrows&{\mathcal G}(Y\times_XY)\end{array}$$
coincide, so the parallel arrows in the lower line also coincide, 
i.e. ${\mathcal G}$ is an ${\mathcal S}$-sheaf.

Assume that there are generically non-finite morphisms in ${\mathcal S}$ 
{\sl with arbitrary targets}. Thus as before, ${\mathcal I}_G$ is a 
particular case of $\Sm_G^{{\mathcal S}}$. Moreover, as restriction of any 
morphism $X\stackrel{f}{\longrightarrow}Y$ to an open dense subset $U$ 
of $X$ factors through $U\stackrel{(f,\phi)}{\longrightarrow}
Y\times{\mathbb A}^1\stackrel{{\rm pr}_Y}{\longrightarrow}Y$, one has 
$\Sm_G^{{\mathcal S}}\subseteq{\mathcal I}_G$. 

1. The categories $\Sm_G^{{\mathcal S}}$ and $\Sm_G$ are abelian, complete, 
cocomplete and have enough injectives. (This is standard.)

2. The section functors $\Hom_{\Sm_G}(\Psi_Y,-):
{\mathcal F}\mapsto{\mathcal F}(Y)$ are exact on $\Sm_G^{{\mathcal S}}$ 
for all smooth $k$-varieties $Y$. As a consequence, quotients of 
${\mathcal S}$-sheaves by their subsheaves coincide with their quotients 
as presheaves: if ${\mathcal F}\in\Sm_G^{{\mathcal S}}$ and ${\mathcal G}$ 
is a subsheaf of ${\mathcal F}$ then 
$({\mathcal F}/{\mathcal G})(Y)={\mathcal F}(Y)/{\mathcal G}(Y)$. 

3. A sheaf is an ${\mathcal S}$-sheaf if and only if all its irreducible 
subquotients are ${\mathcal S}$-sheaves. 

[Proof of the ``only if'' part. As it was shown above, a subsheaf 
${\mathcal G}$ of ${\mathcal F}\in\Sm_G^{{\mathcal S}}$ is an 
${\mathcal S}$-sheaf. By property 2, $({\mathcal F}/{\mathcal G})(Y)=
{\mathcal F}(Y)/{\mathcal G}(Y)$, which implies that the quotient 
${\mathcal F}/{\mathcal G}$ is also an ${\mathcal S}$-sheaf. The ``if'' 
part is shown in Proposition \ref{harakteriz-S-objektov} (in the 
language of representations); cf. also Theorem \ref{irr-subq-hom-inv}.] 

4. The inclusion ${\mathcal I}_G\hookrightarrow\Sm_G$ 
admits a left adjoint ${\mathcal I}$ and a right adjoint. 

Examples of calculation of these adjoint functors are given 
in Propositions \ref{wedge} and \ref{A1-quotient-of-closed}. 

5. The sheaves $C_{k(X)}:={\mathcal I}\Psi_X$ form a system of projective 
generators of ${\mathcal I}_G$. [This follows from 2 and 4.]

({\it Remark.} There are no projective objects in $\Sm_G$.) 

\subsection{Irreducible objects of ${\mathcal I}_G$} \label{Irred-obj} 
{\it Examples.} Let $M$ be a simple effective primitive pure covariant 
motive. Then 
$${\mathbb B}^0(M):Y\mapsto\Hom_{\{\text{{\rm pure $k$-motives}}\}}
(\overline{Y},M)$$ is a well-defined sheaf of finite-dimensional 
${\mathbb Q}$-vector spaces (\cite{repr}). 

A particular case of this example is the sheaf ${\mathcal H}_1^A$, 
corresponding to the motive ``$H_1(A)$'' for any simple abelian 
$k$-variety $A$. 
\begin{proposition}[\cite{repr}] \label{funktor-B} ${\mathbb B}^0$ 
gives rise to a fully faithful functor ${\mathbb B}^{\bullet}$:
$$\{\mbox{{\rm pure $k$-motives}}\}\longrightarrow
\{\mbox{{\rm semisimple sheaves of finite length of 
finite-dimensional graded ${\mathbb Q}$-vector spaces}}\}.$$ 
\end{proposition}
\begin{conjecture}[\cite{repr}] \label{equiv-B} This is an equivalence of 
categories. {\rm (In other words, any irreducible sheaf of finite-dimensional 
${\mathbb Q}$-vector spaces is isomorphic to ${\mathbb B}^0(M)$ 
for a primitive irreducible effective pure motive $M$.)}
\end{conjecture} 

This can be complemented by the following conjecture, which I consider 
as one of the principal problems on ${\mathbb A}^1$-invariant sheaves. 

\begin{conjecture}[\cite{pgl}] \label{irred-diff-for} Any simple 
${\mathbb A}^1$-invariant sheaf can be embedded into the sheaf 
$\SOm^{\bullet}_{|k}:Y\mapsto\Omega^{\bullet}_{k(Y)|k}$. \end{conjecture} 

This conjecture is rather strong: it implies the Bloch's conjecture: 
\begin{conj-corollary}[\cite{pgl}] \label{gip-bloha} 
Suppose that a rational map $f:Y\dasharrow X$ of smooth proper 
$k$-varieties induces an injection 
$\Gamma(X,\Omega^{\bullet}_{X|k})\hookrightarrow
\Gamma(Y,\Omega^{\bullet}_{Y|k})$.\footnote{{\it Example.} Let $r\ge 1$ 
be an integer and $X$ be a smooth proper $k$-variety with 
$\Gamma(X,\Omega^j_{X|k})=0$ for all $r<j\le\dim X$. Let $Y$ be 
a sufficiently general $r$-dimensional plane section of a smooth projective 
variety $X'$ birational to $X$. Then, as all considered invariants 
are birational, the inclusion $Y\hookrightarrow X'$ induces 
an injection $\Gamma(X,\Omega^{\bullet}_{X|k})\hookrightarrow
\Gamma(Y,\Omega^{\bullet}_{Y|k})$.} Then $f$ induces a surjection 
$CH_0(Y)\to CH_0(X)$. 

If $\Gamma(X,\Omega^{\ge 2}_{X|k})=0$ then the Albanese map induces 
an isomorphism $CH_0(X)^0\stackrel{\sim}{\longrightarrow}{\rm Alb}(X)(k)$, 
where $CH_0(X)^0$ is the Chow group of 0-cycles of degree 0 and 
${\rm Alb}(X)$ is the Albanese variety of $X$. 
{\rm (The converse, due to Mumford, is well-known.)} 
In that case $C_{k(X)}=CH_0(X_F)_{{\mathbb Q}}$. 
\end{conj-corollary}
{\it Proof.} Let $C$ be the cokernel of 
$\alpha:CH_0(Y_F)_{{\mathbb Q}}\to CH_0(X_F)_{{\mathbb Q}}$. Then the kernel 
of the homomorphism $\alpha^{\ast}:\Hom_G(CH_0(X_F),
\Omega^{\bullet}_{F|k})\to\Hom_G(CH_0(Y_F),\Omega^{\bullet}_{F|k})$ is 
$\Hom_G(C,\Omega^{\bullet}_{F|k})$. By Proposition \ref{wedge}, the 
homomorphism $\alpha^{\ast}$ coincides with the pull-back under $f^{\ast}
:\Gamma(X,\Omega^{\bullet}_{X|k})\to\Gamma(X,\Omega^{\bullet}_{Y|k})$. 
As the latter is injective, we conclude that 
$\Hom_G(C,\Omega^{\bullet}_{F|k})=0$. If $C\neq 0$ then it 
is cyclic, and thus, admits an simple quotient, and therefore, 
a non-zero morphism to $\Omega^{\bullet}_{F|k}$. 
This contradiction implies that $C=0$. 

As the objects ${\mathbb Q}$ and ${\rm Alb}X(F)_{{\mathbb Q}}$ of 
${\mathcal I}_G$ are projective (\cite[\S6.2]{repr}), the natural 
surjections $\deg:C_{k(X)}\to{\mathbb Q}$ and ${\rm Alb}_F:\ker\deg\to
{\rm Alb}X(F)_{{\mathbb Q}}$ are split, so the cyclic $G$-module 
$C_{k(X)}$ is isomorphic to a direct sum of type ${\mathbb Q}\oplus
{\rm Alb}X(F)_{{\mathbb Q}}\oplus\ker{\rm Alb}_F$. Thus, 
$\Hom_G(C_{k(X)},\Omega^{\bullet}_{F|k})\cong\Hom_G({\mathbb Q}\oplus
{\rm Alb}X(F),\Omega^{\bullet}_{F|k})\oplus\Hom_G(\ker{\rm Alb}_F,
\Omega^{\bullet}_{F|k})$. By Proposition \ref{wedge}, 
$\Hom_G(C_{k(X)},\Omega^{\bullet}_{F|k})=\Gamma(X,\Omega^{\bullet}_{X|k})$ 
and $\Hom_G({\mathbb Q}\oplus{\rm Alb}X(F),\Omega^{\bullet}_{F|k})=
\Gamma(X,\Omega^{\le 1}_{X|k})$. If $\Gamma(X,\Omega^{\ge 2}_{X|k})=0$ this 
means that $\Hom_G(C_{k(X)},\Omega^{\bullet}_{F|k})=\Hom_G({\mathbb Q}\oplus
{\rm Alb}X(F),\Omega^{\bullet}_{F|k})$. Therefore, 
the $G$-module $\ker{\rm Alb}_F$ should be zero, as otherwise it 
is cyclic, thus admits a non-zero simple quotient, and (by Conjecture 
\ref{irred-diff-for}) a non-zero morphism to $\Omega^{\bullet}_{F|k}$. 
It remains to take the $G$-invariants of $\ker\deg
\stackrel{\sim}{\longrightarrow}CH_0(X_F)^0_{{\mathbb Q}}
\stackrel{\sim}{\longrightarrow}{\rm Alb}X(F)_{{\mathbb Q}}$; 
the torsion is controlled by Roitman's theorem. \qed 

\vspace{4mm}

Also this would imply that any irreducible ${\mathbb A}^1$-invariant 
sheaf is a sheaf of finite-dimensional vector spaces. 

{\it Example.} Let ${\mathcal F}$ be a simple ${\mathbb A}^1$-invariant 
sheaf and suppose that it is of level 1, i.e. it is non-constant and 
${\mathcal F}(Y)\neq 0$ for a curve $Y$, cf. also 
p.\pageref{coniveau-filtr}. Then, by \cite[Corollary 6.22]{repr}, 
${\mathcal F}\cong{\mathcal H}^A_1$ for a simple abelian variety $A$. 
Now any non-zero $\eta\in\Gamma(A,\Omega^1_{A|k})$ gives an embedding 
${\mathcal F}\hookrightarrow\Omega^1_{|k}$ by $[x:{\mathcal O}(U)\to k(Y)]
\mapsto x(\eta)\in\Omega^1_{Y|k}(Y)$ ($U\subset A$ is an affine open subset). 

\begin{proposition} \label{bir-A1-presheaf} A dominant presheaf 
${\mathcal F}$ is a sheaf if and only if the following three conditions hold: 
{\rm (i)} the sequence ${\mathcal F}(X)\to{\mathcal F}(X\times{\mathbb A}^1)
\rightrightarrows{\mathcal F}(X\times{\mathbb A}^2)$ is exact for any 
smooth $k$-variety $X$,\footnote{E.g., any ${\mathbb A}^1$-invariant presheaf 
${\mathcal F}$ satisfies the condition (i).} {\rm (ii)} ${\mathcal F}$ is 
birationally invariant, {\rm (iii)} it has the Galois descent property, 
i.e. ${\mathcal F}(X)={\mathcal F}(Y)^{{\rm Aut}(Y|X)}$ 
for any Galois covering $Y\to X$. \end{proposition}
{\it Proof.} The conditions (i)--(iii) are particular cases of the equalizer 
diagram (\ref{opr-puchka}) for coverings by (i) projections 
$X\times{\mathbb A}^s\to X$, (ii) open dense $U\subset X$, (iii) \'etale 
Galois covers $Y\to X$, respectively. Galois descent property for any 
sheaf is clear, since \'etale morphisms with dense images are covering 
and $U\times_XU=\coprod_{g\in{\rm Aut}(Y|X)}U_g$ for a Zariski open 
${\rm Aut}(Y|X)$-invariant $U\subset Y$, where $U_g\cong U$ is 
the image of the embedding $(id_U,g):U\hookrightarrow U\times_XU$. 

Conversely, it is clear that any Galois-separable presheaf ${\mathcal F}$ 
satisfying (i) and (ii) is separable: if $Y\to X$ is a cover, i.e. 
a smooth dominant morphism, then for any sufficiently general dominant map 
$\varphi:Y\dasharrow{\mathbb A}^{\delta}$ (where $\delta=\dim Y-\dim X$) 
we can choose a dominant \'etale morphism $\widetilde{Y}\to Y$ so that the 
composition $\widetilde{Y}\to Y\dasharrow X\times{\mathbb A}^{\delta}$ 
is Galois with the group denoted by $H$, and therefore, the composition 
$$\begin{array}{rccc}{\mathcal F}(X)
\stackrel{\text{(i)}}{\hookrightarrow}
{\mathcal F}(X\times{\mathbb A}^1) 
\stackrel{\text{(i)}}{\hookrightarrow}
\dots\stackrel{\text{(i)}}{\hookrightarrow}&
{\mathcal F}(X\times{\mathbb A}^{\delta})&\longrightarrow&
{\mathcal F}(Y)\\ &\downarrow\lefteqn{\text{injective}}&&\downarrow\\
&{\mathcal F}(\widetilde{Y})^H&\hookrightarrow&
{\mathcal F}(\widetilde{Y})\end{array}$$ is injective. 
Then in the commutative diagram 
\begin{equation}\label{red-to-aff}
\begin{array}{ccccc}{\mathcal F}(X)&\to&{\mathcal F}(\widetilde{Y})
&\rightrightarrows&{\mathcal F}(\widetilde{Y}\times_X\widetilde{Y})\\
\|&&\uparrow&&\uparrow\\
{\mathcal F}(X)&\to&{\mathcal F}(Y)&\rightrightarrows&
{\mathcal F}(Y\times_XY)\\
\|&&\uparrow&&\uparrow\\
{\mathcal F}(X)&\to&{\mathcal F}(X\times{\mathbb A}^{\delta})&
\rightrightarrows&
{\mathcal F}(X\times{\mathbb A}^{\delta}\times{\mathbb A}^{\delta})
\end{array}\end{equation} all arrows are injective, so it suffices 
to show the exactness of the upper row. 

Let $f$ be an element of ${\mathcal F}(\widetilde{Y})$. The image of 
$f$ in ${\mathcal F}(\widetilde{Y}\times_X\widetilde{Y})$ under the 
projection to the first factor is fixed by $\{1\}\times H$; the image 
of $f$ in ${\mathcal F}(\widetilde{Y}\times_X\widetilde{Y})$ under 
the projection to the second factor is fixed by $H\times\{1\}$. Now 
if $f$ is an element of the equalizer of ${\mathcal F}(\widetilde{Y})
\rightrightarrows{\mathcal F}(\widetilde{Y}\times_X\widetilde{Y})$ then 
the two images coincide, so they are fixed by the group $H\times H$. 
The injectivity of both parallel arrows in the upper row of the diagram 
(\ref{red-to-aff}) implies that $f\in{\mathcal F}(\widetilde{Y})^H$. 
By (iii) and the injectivity of the vertical arrow, $f$ comes from the 
equalizer of the bottom row of the diagram (\ref{red-to-aff}). 
Finally, the bottom row of the diagram (\ref{red-to-aff}) is exact by 
Lemma \ref{red-to-A1}, and thus, $f$ comes from ${\mathcal F}(X)$. \qed

\begin{lemma} \label{red-to-A1} Let ${\mathcal V}$ be a category of 
schemes such that for any $X\in{\mathcal V}$: {\rm (i)} the projection 
$X\times{\mathbb A}^1\to X$ is a morphism in ${\mathcal V}$, 
{\rm (ii)} any linear automorphism of any affine space ${\mathbb A}$ 
induces an automorphism of $X\times{\mathbb A}$ in ${\mathcal V}$. 
Let ${\mathcal F}$ be a presheaf on this category such that 
the sequence ${\mathcal F}(X)\to{\mathcal F}(X\times{\mathbb A}^1)
\rightrightarrows{\mathcal F}(X\times{\mathbb A}^2)$ is exact 
for any $X\in{\mathcal V}$. Then the sequence ${\mathcal F}(X)
\to{\mathcal F}(X\times{\mathbb A}^s)\rightrightarrows
{\mathcal F}(X\times{\mathbb A}^{2s})$ is exact for any 
$X\in{\mathcal V}$ and any $s\ge 1$. \end{lemma} 
{\it Proof.} We proceed by induction on $s$, the case $s=1$ being trivial. 

Denote by $\mathrm{pr}_1,\mathrm{pr}_2:X\times{\mathbb A}^{2s}
\rightrightarrows X\times{\mathbb A}^s$ the two projections. 
 
For any $f\in{\mathcal F}(X\times{\mathbb A}^s)$ the element 
$\mathrm{pr}_1^{\ast}f$ is fixed by $\Phi^{\ast}\in\End
{\mathcal F}(X\times{\mathbb A}^s\times{\mathbb A}^s)$ 
for any linear automorphism $\Phi(u,v)=(u,\varphi(u,v))$ of 
${\mathbb A}^s\times{\mathbb A}^s$. 
Similarly, $\mathrm{pr}_2^{\ast}f$ is fixed by $\Psi^{\ast}\in\End
{\mathcal F}(X\times{\mathbb A}^s\times{\mathbb A}^s)$ 
for any linear automorphism $\Psi(u,v)=(\psi(u,v),v)$. 

Let now $f\in{\mathcal F}(X\times{\mathbb A}^s)$ be in the 
equalizer of $\mathrm{pr}_1^{\ast}$ and $\mathrm{pr}_2^{\ast}$. 
Then $\mathrm{pr}_1^{\ast}f=\mathrm{pr}_2^{\ast}f$ is fixed by the 
group, generated by $\Phi^{\ast}$ and $\Psi^{\ast}$ as above. 
Clearly, such automorphisms $\Phi$ and $\Psi$ generate the group 
consisting of all linear automorphisms $\alpha$. Then 
$\mathrm{pr}_1^{\ast}f=\alpha^{\ast}\mathrm{pr}_1^{\ast}f$. 

Applying the induction assumption in the case where $\alpha$ is 
identical on one of the first $s$ coordinates and interchanges 
$i$-th and $(s+i)$-th for other $1\le i\le s$, we get that 
$f$ belongs to the image of ${\mathcal F}(X\times{\mathbb A}^1)\to
{\mathcal F}(X\times{\mathbb A}^s)$ under morphism induced 
by the projection ${\mathbb A}^s\to{\mathbb A}^1$ to one 
of the copies of ${\mathbb A}^1$. Then the case $s=1$ implies 
that $f$ comes from  ${\mathcal F}(X)$.  \qed 

\vspace{4mm}

{\sc Examples.} 1.  A stable birationally invariant dominant presheaf 
with the Galois descent is a sheaf. 

2. Example of a birationally invariant presheaf ${\mathcal F}$ with 
the Galois descent property which is not a sheaf. Let ${\mathcal G}$ 
be a dominant sheaf and $I\subsetneq\{0,1,2,\dots\}$ be a non-empty 
(finite or infinite) interval. Assume that ${\mathcal G}(X)\neq 0$ 
for some $X$ with $\dim X\not\in I$. Then the presheaf ${\mathcal F}:
U\mapsto\left\{\begin{array}{ll}{\mathcal G}(U)&\text{if $\dim U\in I$}\\ 
0&\text{if $\dim U\not\in I$}\end{array}\right.$ (with the restriction maps 
of ${\mathcal G}$, whenever possible, otherwise zero) is birationally 
invariant and has the Galois descent property, but it is not a sheaf. The 
sheafification of ${\mathcal F}$ is ${\mathcal G}$ if $I$ is infinite 
and 0 otherwise. 

\vspace{4mm}

Now, what are the projective generators of ${\mathcal I}_G$ from 
\S\ref{Properties-of}, Property 5?
\begin{conjecture} \label{C=CH} For any smooth proper $k$-variety $X$,
the sheaf $C_{k(X)}$ coincides with $Y\mapsto CH_0(X_{k(Y)})_{{\mathbb Q}}$.
\end{conjecture}

{\it Remarks.} 1. This is known, e.g., if $X$ is a curve, cf. 
\cite[Cor.6.21]{repr} and Proposition \ref{pic-as-adj} for a stronger 
statement. Conjecture \ref{C=CH} would imply that ${\mathcal I}_G$ is a 
tensor category under the operation $({\mathcal F},{\mathcal G})\mapsto
{\mathcal I}({\mathcal F}\otimes_{\Sm_G}{\mathcal G})=:{\mathcal F}
\otimes_{{\mathcal I}}{\mathcal G}$, where $\otimes_{\Sm_G}$
denotes the sheafification of the tensor product presheaf.
Moreover, the ``K\"{u}nneth formula'' holds:
$C_{k(X)}\otimes_{{\mathcal I}}C_{k(Y)}=C_{k(X\times_kY)}$.

2. It is shown in \cite[Proposition 6.17]{repr} that, roughly speaking, 
$C_{k(X)}$ is the quotient of generic 0-cycles on $X$ by those divisors 
of rational functions on generic curves on $X$ which are generic, 
and thus, Conjecture \ref{C=CH} should be considered as a moving lemma. 

3. Conjecture \ref{C=CH} and the motivic conjectures imply 
conjectures \ref{equiv-B} and \ref{irred-diff-for}. 

\section{Alternative descriptions of ${\mathbb A}^1$-invariant sheaves}
\label{Alternative-descr} 
Now I want to introduce the language of representations and to use 
it to explain some results and conjectures of \S\ref{pre-sheaves}, 
especially Conjecture \ref{C=CH}. 

\subsection{Smooth representations and non-degenerate modules over 
algebras of measures} 
For any totally disconnected Hausdorff group\footnote{cf. 
\cite[Appendix A]{glatt}} $H$ an $H$-set (group, etc.) is called 
{\sl smooth} if the stabilizers are open. 

Any smooth representation $W$ of $H$ over $E$ can be considered as 
a module over the associative algebra ${\mathbb D}_E(H):=\prlim_UE[H/U]$ 
of the ``oscillating'' measures on $H$ (for which all open subgroups and 
their translates are measurable): ${\mathbb D}_E(H)\times W\to W$ is 
defined by $(\alpha,w)\mapsto\beta w$ for any $\beta\in E[H]$ with the 
same image $E[H]$ as $\alpha$, where $w\in W^U$ for some open subgroup 
$U$ of $H$. 

Passing to the inverse limit, we get the algebra structure 
on ${\mathbb D}_E(H)$ from ${\mathbb D}_E(H)\times E[H/U]\to E[H/U]$. 

\label{nul-annul-sod-vse} If the annihilator of $W\in\Sm_H(E)$ in 
${\mathbb D}_E(H)$ vanishes then the restriction of $W$ to any compact 
subgroup $U$ contains each smooth irreducible representation of $U$. 
(Otherwise, if $W$ does not contain a smooth irreducible representation 
$\rho$ of $U$ then the natural projector in ${\mathbb D}_E(H)$ to the 
$\rho$-isotypical part would annihilate $W$.) 

\subsection{A representation theoretic setting for (${\mathbb A}^1$-invariant) 
sheaves} \label{gomotopich-invar-prdst} In this section, for a group $H$ 
as in \S\ref{nul-annul-sod-vse} and a collection $S$ of pairs of its 
subgroups, we study the category $\Sm_H^S(E)$ of smooth $E$-representations 
$W$ of $H$, satisfying $W^{U_1}=W^{U_2}$ for all $(U_1,U_2)\in S$. 

Theorem \ref{sheafif-i-equiv} explains the consistence 
of this notation with that of \S\ref{pre-sheaves}. 

Collections $S$ and $S'$ are called {\sl equivalent} 
if they define the same subcategory of $\Sm_H:=\Sm_H^{\emptyset}$. 

For any subgroup $U\subset H$ the functor $H^0(U,-)$ on the category of 
smooth $H$-sets (or modules, etc.) coincides with $\indlim H^0(V,-)$, where 
the limit is taken over the open 
subgroups $V$ of $H$ containing $U$. Therefore, one can assume that the 
subgroups $U_1,U_2$ are intersections of open ones, and in particular, that 
they are closed. Further, as $W^{U_1}\cap W^{U_2}=W^{\langle U_1,U_2\rangle}$ 
for any $H$-module $W$ and $U_1,U_2\subseteq\langle U_1,U_2\rangle$, one can 
assume that the pairs $(U_1,U_2)\in S$ are ordered: $U_1\subset U_2$. 

\begin{lemma} \label{podfaktory-dlq-komp} 
Assume that for any pair $(U_1\subset U_2)\in S$ the functor $H^0(U_1,-)$ 
is exact on $\Sm_H$. Then the category $\Sm_H^S(E)$ is stable under passing 
to the subquotients in $\Sm_H(E)$, and in particular, it is abelian. 
The inclusion functor $\Sm_H^S(E)\hookrightarrow\Sm_H(E)$ 
admits a left adjoint\footnote{The diagrams $\begin{array}{ccc}\Sm_H(E)&
\stackrel{{\mathcal I}_S}{\longrightarrow}&\Sm^S_H(E)\\ 
\otimes_EE'\downarrow\phantom{\otimes_EE'}&&
\phantom{\otimes_EE'}\downarrow\otimes_EE'\\
\Sm_H(E')&\stackrel{{\mathcal I}_S}{\longrightarrow}&\Sm^S_H(E')\end{array}
\quad\mbox{and}\quad\begin{array}{ccc}\Sm_H(E)&
\stackrel{{\mathcal I}_S}{\longrightarrow}&\Sm^S_H(E)\\ 
{\rm for}\uparrow\phantom{{\rm for}}&&\phantom{{\rm for}}\uparrow{\rm for}\\
\Sm_H(E')&\stackrel{{\mathcal I}_S}{\longrightarrow}
&\Sm^S_H(E')\end{array}$ are commutative for any field extension $E'|E$, 
so omitting $E$ from the notation does not lead to a confusion.} 
$W\longmapsto{\mathcal I}_SW$. \end{lemma}
{\it Proof.} If a sequence $0\to W_1\to W\to W_2\to 0$ in $\Sm_H$ 
is exact then the sequences $0\to W_1^{U_1}\to W^{U_1}\to W_2^{U_1}\to 0$ 
and $0\to W_1^{U_2}\to W^{U_2}\to W_2^{U_2}$ are also exact. If 
$W\in\Sm_H^S$, i.e. $W^{U_1}=W^{U_2}$, then 
$W_1^{U_1}=W_1\cap W^{U_1}=W_1\cap W^{U_2}=W_1^{U_2}$ and 
$W^{U_2}\to W_2^{U_2}$ is surjective (since $W^{U_2}=W^{U_1}\to W_2^{U_1}$ 
is surjective and factors through $W_2^{U_2}\subseteq W_2^{U_1}$).
This means that $\Sm_H^S$ is stable under taking subquotients in $\Sm_H$. 

The existence of the functor ${\mathcal I}_S$ can be deduced from the 
special adjoint functor theorem, cf. \cite[\S5.8]{maclane}. However, we 
construct it ``explicitly'', which enables us to relate the generators 
of the category ${\mathcal I}_G$ to the Chow groups of 0-cycles. 

Let $W'\in\Sm^S_H$. Any $H$-homomorphism 
$W\stackrel{\alpha}{\longrightarrow}W'$ factors through the object 
$\alpha(W)$ of $\Sm^S_H$. We may, therefore, assume that 
$\alpha$ is surjective. Let $(U_1\subset U_2)\in S$. As the functor 
$H^0(U_1,-)$ is exact on $\Sm_H$, the morphism $\alpha$ 
induces a surjection $W^{U_1}\longrightarrow(W')^{U_1}$. As 
$(W')^{U_2}=(W')^{U_1}$, the subgroup $U_2$ acts on $(W')^{U_1}$ trivially, 
and therefore, the subrepresentation 
$W_{U_1\subset U_2}=\langle\sigma w-w~|~\sigma\in U_2,~w\in W^{U_1}\rangle_H$ 
of $H$ is contained in the kernel of $\alpha$. It follows that $\alpha$ 
factors through ${\mathcal I}_SW:=W/\sum_{(U_1\subset U_2)\in S}
W_{U_1\subset U_2}$. 

The representation ${\mathcal I}_SW$ of $H$ is smooth, so the map 
$W^{U_1}\longrightarrow({\mathcal I}_SW)^{U_1}$, induced by the projection, 
is surjective, and therefore, any element $\overline{w}\in
({\mathcal I}_SW)^{U_1}$ can be lifted to an element $w\in W^{U_1}$. 
Then $\sigma\overline{w}-\overline{w}$ coincides with the projection 
of the element $\sigma w-w$ for any $\sigma\in U_2$. Notice that 
$\sigma w-w\in W_{U_1\subset U_2}$, so its projection is zero, and 
therefore, $\sigma\overline{w}=\overline{w}$ for any $\sigma\in U_2$. 
As $({\mathcal I}_SW)^{U_2}\subseteq({\mathcal I}_SW)^{U_1}$, this means 
that $({\mathcal I}_SW)^{U_2}=({\mathcal I}_SW)^{U_1}$, and thus, 
${\mathcal I}_SW\in\Sm^S_H$. 

One has $\Hom_{\Sm^S_H}({\mathcal I}_SW,W')=\Hom_{\Sm_H}(W,W')$ for any 
$W\in\Sm_H$ and $W'\in\Sm^S_H$, i.e. the functor ${\mathcal I}_S$ is left 
adjoint to the inclusion functor $\Sm^S_H\hookrightarrow\Sm_H$. \qed 

\vspace{4mm}

{\it Remark.} The functor ${\mathcal I}_S$ generalizes the coinvariants, 
since ${\mathcal I}_S=H_0(H,-)$ if $S=\{(\{1\}\subset H)\}$. 

{\sc Examples.} 1. The functor $H^0(U_1,-)$ is exact on $\Sm_H$ if, 
e.g., the subgroup $U_1$ is compact. 

2. Suppose that $H$ is the automorphism group of an algebraically closed 
field extension $F|k$ of countable transcendence degree and $U_1$ is the 
subgroup of automorphisms of $F$ over a fixed subextension of $k$ in $F$ 
of infinite transcendence degree. Though $U_1$ need not be 
compact, the functor $H^0(U_1,-)$ is exact on $\Sm_H$. 

\begin{proposition} \label{harakteriz-S-objektov} 
Let $H$ be a totally disconnected group and $S$ be such a collection of 
pairs of its subgroups $(U_1\subset U_2)$ that 
\begin{enumerate}\item \label{cond-1} 
for any pair $(U_1\subset U_2)\in S$ there exists 
an element $\sigma\in U_2$ such that 
{\rm (i)} $(U_1\cap\sigma U_1\sigma^{-1}\subset U_1)\in S$; {\rm (ii)} 
$U_1$ and $\sigma U_1\sigma^{-1}$ generate $U_2$, at least topologically. 
\item there exists an equivalent collection of pairs of its subgroups 
$(U_1\subset U_2)$, where all $U_1$ are compact. \end{enumerate}
Then an object of $\Sm_H(E)$ belongs to $\Sm_H^S(E)$ if 
and only if all its irreducible subquotients are in $\Sm_H^S(E)$. 
In particular, $\Sm_H^S(E)$ is a Serre subcategory of 
$\Sm_H(E)$. \end{proposition} 
{\it Proof.} Suppose that $W\not\in\Sm_H^S$, whereas all its 
irreducible subquotients are in $\Sm_H^S$. Then $W^{U_1}\neq 
W^{U_2}$ for some pair $(U_1\subset U_2)\in S$, that is there exist a 
vector $v\in W^{U_1}\smallsetminus W^{U_2}$. Choose an element 
$\sigma\in U_2$ as in condition (\ref{cond-1}) of the statement for the 
pair $(U_1\subset U_2)\in S$. Then $\sigma v-v=:u\neq 0$, since $U_1$ and 
$\sigma$ generate a dense subgroup in $U_2$. 

One may replace $W$ by its quotient by a maximal subrepresentation not 
containing $u$. Then the subrepresentation $\langle u\rangle$, generated 
by $u$, becomes irreducible, and thus, an object of $\Sm_H^S$. 

By definition, $u\in W^{U_1}+W^{\sigma U_1\sigma^{-1}}\subseteq 
W^{U_1\cap\sigma U_1\sigma^{-1}}$. As $\langle u\rangle\in\Sm_H^S$ 
and $(U_1\cap\sigma U_1\sigma^{-1}\subset U_1)\in S$, we conclude that 
$u\in W^{U_1}$. This implies that $\sigma v\in W^{U_1}$. On the other hand, 
$\sigma v\in W^{\sigma U_1\sigma^{-1}}$, so $\sigma v\in W^{U_1}\cap 
W^{\sigma U_1\sigma^{-1}}$. The latter vector space coincides with $W^{U_2}$, 
and thus, $v\in W^{U_2}$, contradicting our assumptions. 

The converse follows from Lemma \ref{podfaktory-dlq-komp}. \qed 

\subsection{More notations and compatibility of notations of 
\S\ref{gomotopich-invar-prdst} and \S\ref{pre-sheaves}: 
the sheafification and smooth representations} \label{obozn-predst} 
From now on we fix the following notations: $F|k$ is an algebraically 
closed field extension of countably infinite transcendence degree, 
and $G=G_{F|k}$ is the automorphism group of the extension $F|k$. 

Consider connected smooth $k$-varieties $U$ endowed with a generic 
$F$-point, i.e., with a $k$-field embedding $k(U)\overk F$. 
For any presheaf ${\mathcal F}$ on ${\mathcal S}m_k'$ we can form the 
direct limit ${\mathcal F}(F):=\indlim{\mathcal F}(U)$ over such $U$. 
The group $G={\rm Aut}(F|k)$ acts naturally on ${\mathcal F}(F)$. 
\begin{theorem}[\cite{glatt}]\label{sheafif-i-equiv} \begin{itemize}\item 
${\mathcal F}\mapsto{\mathcal F}(F)$ gives an equivalence of the categories 
$\Sm_G^{{\mathcal S}}(E)$ (of \S\ref{pre-sheaves}) and $\Sm_G^S(E)$ 
(of \S\ref{gomotopich-invar-prdst}), where $S$ is the collection 
of pairs $G_{F|k(X)}\subseteq G_{F|k(Y)}$ for all morphisms 
$(X\to Y)\in{\mathcal S}$. 
\item For any presheaf ${\mathcal F}$, the sheaf corresponding to 
${\mathcal F}(F)$ is the sheafification of ${\mathcal F}$. 
\end{itemize} \end{theorem}

\subsection{An example: birational invariants constant on the projective 
spaces} Let $S$ consist of a single pair $K\subset G$ such that $K$ is a 
`maximal' compact subgroup, i.e., any compact subgroup is conjugate to a 
subgroup of $K$. Then $S$ is equivalent to the collection consisting of a 
single pair $K'\subset G$, where $K'$ is the pointwise stabilizer of some 
transcendence base of $F|k$, and also to the collection $S'$ of pairs 
$U\subset G$ such that $U$ is the pointwise stabilizer of a finite subset 
of a fixed transcendence base of $F|k$. The collection $S'$ satisfies 
the assumptions of Proposition \ref{harakteriz-S-objektov}. 

\begin{lemma} \label{nepriv-tenz-proizv} Let $E'|E$ be an extension 
of fields, $H$ be a group and $(\rho,W_2)$ be an irreducible 
$E'$-representation of $H$. Let $W_1$ be an $E$-representation of $H$, 
absolutely irredicible even in restriction to $\ker\rho$.\footnote{i.e., 
irredicible and with ${\rm End}_{E[\ker\rho]}(W_1)=E$: otherwise, if 
${\rm End}_{E[\ker\rho]}(W_1)\neq E$ and $E'|E$ is a non-trivial 
field extension in the division $E$-algebra ${\rm End}_{E[\ker\rho]}(W_1)$ 
then the action of $E'$ on $W_1$ gives a non-injective surjection of 
$E'$-representations $W_1\otimes_EE'\longrightarrow W_1$.} Then the 
$E'$-representation $W_1\otimes_EW_2$ of $H$ is irredicible. \end{lemma} 
{\it Proof.} Let $\xi\in W_1\otimes_EW_2$ be a non-zero vector. 
It suffices to check that the $E'[H]$-span of $\xi$ contains 
$W_1\otimes v$ for any $v\in W_2$. 
Any non-zero $E[\ker\rho]$-submodule in $W_1^m$ is isomorphic to $W_1^{m'}$ 
for some $1\le m'\le m$, and therefore, the $E[\ker\rho]$-submodule 
in $W_1\otimes_EW_2$ spanned by $\xi$ (which is in fact a submodule 
in $\oplus_{i=1}^mW_1\otimes v_i\cong W_1^m$ for some $m\ge 1$ and 
$E$-linearly independent $v_1,\dots,v_m$) contains 
a $E[\ker\rho]$-submodule $W_1'$ isomorphic to $W_1$. 
As the endomorphisms of the $E[\ker\rho]$-module $W_1$ are scalar, 
there exists a non-zero $m$-tuple $(a_1,\dots,a_m)\in E^n$ such that 
$W_1'=\{a_1w\otimes v_1+\dots+a_mw\otimes v_m~|~w\in W_1\}$. 
In other words, $W_1'=W_1\otimes v'$, where $v':=a_1v_1+\dots+ a_mv_m$ 
is a non-zero vector in $W_2$. 

As any vector $v$ of $W_2$ is an $E'$-linear combination of several 
elements in the $H$-orbit of $v'$, we may assume that $v=hv'$ for some 
$h\in H$. Then $u\otimes v=h(h^{-1}u\otimes v')$ for any $u\in W_1$. \qed 

\begin{lemma} \label{max-v-tenz-proizv} Let $E'|E$ be an extension 
of fields, $H$ be a group and $(\rho,W_2)$ be an irreducible 
$E'$-representation of $H$. Let $W_1$ be an $E$-representation of $H$ 
such that {\rm (i)} the sum $\Sigma$ of all proper $E$-subrepresentations 
of $\ker\rho$ in $W_1$ is proper\footnote{$\Sigma$ is $H$-invariant: 
as $\ker\rho$ is a normal subgroup of $H$, the group $H$ permutes the 
$\ker\rho$-submodules in $W_1$, while $\Sigma$ is the maximal proper 
$\ker\rho$-submodule in $W_1$.} and {\rm (ii)} $W_1/\Sigma$ is 
absolutely irredicible in restriction to $\ker\rho$ and its restriction 
to the pointwise stabilizer $\Xi$ of $\Sigma$ in $\ker\rho$ is non-trivial. 
Then any proper $E'$-subrepresentation of $H$ in $W_1\otimes_EW_2$ is 
contained in $\Sigma\otimes_EW_2$. \end{lemma} 
{\it Proof.} Let $\xi\in W_1\otimes_EW_2$ be a vector, which is not in 
$\Sigma\otimes_EW_2$. It suffices to check that the $E'[H]$-span $V$ of 
$\xi$ contains $W_1\otimes v$ for some non-zero $v\in W_2$, as then $V$ 
coincides with $W_1\otimes_EW_2$: any vector of $W_2$ is an $E'$-linear 
combination of several elements in the $H$-orbit of $v$ and 
$W_1\otimes hv=h(W_1\otimes v)$ for any $h\in H$. 

It follows from Lemma \ref{nepriv-tenz-proizv} that $V$ is surjective over 
$(W_1/\Sigma)\otimes_EW_2$. In particular, $V$ contains an element of type 
$\sum_{i=1}^ma_i\otimes b_i$ for some $a_1\in W_1\smallsetminus\Sigma$, whose 
projection to $W_1/\Sigma$ is not fixed by $\Xi$, for some $a_2,\dots,a_m\in
\Sigma$ and for some $E'$-linearly independent $b_1,\dots,b_m\in W_2$. 
Then there exists an element $h\in\Xi$ such that $ha_1-a_1\in 
W_1\smallsetminus\Sigma$, and therefore, $V$ contains an element of type 
$\sum_{i=1}^ma\otimes b_1$ for some $a\in W_1\smallsetminus\Sigma$. \qed 

\begin{proposition} \label{I-S-polulin} Let $W\in\Sm_G(E)$ be an object. 
For any open subgroup $U$ of $G$, denote by $W_{(U)}$ the sum of all proper 
subrepresentations of $U$ in $W$; and by $\Xi_U$ the pointwise stabilizer of 
$W_{(U)}$ in $U$. Suppose that for any open subgroup $U$ of $G$: {\rm (i)} 
the $E$-representation $W/W_{(U)}$ of $U$ is absolutely irreducible and 
non-trivial in restriction to $\Xi_U$\footnote{In particular, $W$ is absolutely 
indecomposable. Any non-zero quotient of $A(F)$ for an absolutely simple 
algebraic $k$-group $A$ is an example of such $W$. (Indeed, any open subgroup 
$U\subset G$ contains $G_{F|L}$ for a finitely generated $L$ in $F|k$, so 
any $t\in A(F)\smallsetminus A(\overline{L})$ is a cyclic vector of $A(F)$, 
considered as $U$-module. Here $\overline{L}$ is the algebraic closure of 
$L$ in $F$. If the transcendence degree of $L|k$ is minimal then, by 
\cite{max}, $\overline{L}$ is $U$-invariant, so $A(F)_{(U)}=A(\overline{L})$. 
\qed)} and {\rm (ii)} any irreducible smooth representation of $K$ can be 
embedded into $W$ so that its image does not meet $W_{(U)}$. Then 
${\mathcal I}_S$ annihilates any quotient of $W\otimes_EV$ for any 
$V\in\Sm_G(E)$. \end{proposition} 
{\it Proof.} It suffices to check the vanishing of 
${\mathcal I}_S(W\otimes_EV)$. Extending the coefficients if needed, 
we may assume that $E$ is big enough (i.e., algebraically closed and 
$\# E>\# k$), so that any smooth irreducible $E$-representation of 
any open subgroup of $G$ is absolutely irreducible.\footnote{Schur's 
lemma=\cite[Claim 2.11]{bz}: {\it Let $H$ be a totally disconnected 
group and $E$ be a field of cardinality greater than the cardinality of 
$H/U$ for any open subgroup $U$ of $H$. Then the endomorphisms of the 
smooth irreducible $\overline{E}$-representation of $H$ are scalar.}} 

The vanishing holds if the $G$-module $W\otimes_EV$ is spanned by the 
elements $g\xi-\xi$ for all $\xi\in(W\otimes_EV)^K$ and all $g\in G$. 
Equivalently, as the restriction of $V$ to $K$ is semisimple, the $G$-span 
of such elements $g\xi-\xi$ contains $W\otimes_E\rho$ for any irreducible 
$E$-subrepresentation $\rho$ of $K$ in $V$. 
By (ii), $W$ contains a $E$-subrepresentation of $K$ which is 
(a) dual to $\rho$ and (b) outside of $W_{(U)}$, where 
$U\subset G$ is the pointwise stabilizer of $\rho$. 
Then there is an element $\xi\in(W\otimes_E\rho)^K$, which is not in 
$W_{(U)}\otimes_E\rho$. 

As the $\Xi_U$-module $W/W_{(U)}$ is non-trivial, there exists an element 
$u\in\Xi_U$ such that $\eta:=u\xi-\xi$ is not in $W_{(U)}\otimes_E\rho$. 
Denote by $\widetilde{U}$ the subgroup in $G$ generated by $U$ and $K$. 
Then $\widetilde{U}$ contains $U$ as a normal subgroup of finite index; 
$\widetilde{U}$ acts on $W_{(U)}$; $\rho$ can be viewed as a representation 
of $\widetilde{U}$ via the identification $\widetilde{U}/U=K/U\cap K$. 
By Lemma \ref{max-v-tenz-proizv} (with $H=\widetilde{U}$), the element 
$\eta$ generates the $E[\widetilde{U}]$-module $W\otimes_E\rho$. \qed 

\begin{lemma}[A source of representations of $G$ containing all 
irreducible smooth representations of $K$] \label{primery-soderzh-vse}
If a subrepresentation $W$ of $G$ in 
$\bigotimes^{\bullet}_F\Omega^1_{F|k}$ does not contain regular 
forms,\footnote{Examples of such $W$ are subrepresentations of 
${\rm Sym}^2_F\Omega^1_{F|k}$, of $\Omega^{\bullet}_{F|k,\text{exact}}$, 
or of the image in $\Omega^j_{F|k}$ of $\wedge^j\Omega^1_{F|k,\log}$, 
where $d\log:F^{\times}/k^{\times}\stackrel{\sim}{\longrightarrow}
\Omega^1_{F|k,\log}$, for any $j\ge 1$. It follows directly 
from Hilbert's Satz 90 that the representation $F$ (and therefore, the 
irreducible representation $d:F/k\stackrel{\sim}{\longrightarrow}
\Omega^1_{F|k,\text{exact}}$) of $G$ contains all irreducible smooth 
(and thus, finite-dimensional) representations of $K$.} i.e., 
forms from $\Gamma(X,\Omega^{\bullet}_{X|k})$ for a smooth proper 
$k$-variety $X$ with $k(X)\subset F$, then $W$ contains each 
irreducible smooth representation $\rho$ of $K$. \end{lemma} 

As mentioned in \S\ref{nul-annul-sod-vse}, if no non-zero element 
of ${\mathbb D}_E(G)$ annihilates a smooth 
representation $W$ then $W$ contains all irreducible smooth 
representations of $K$. The vanishing of the annihilators of 
$F/k$ and $F^{\times}/k^{\times}$ is shown in \cite[Prop.4.2]{repr}. 
Assume for simplicity that $F^K|k$ is purely transcendental. 

{\it Proof.} Let $p_{\rho}$ be the central projector 
in the group algebra of $Q=K/\ker\rho$ onto the $\rho$-isotypical part. 
As explained in \cite[Prop.7.6]{pgl}, $W$ contains a non-zero element 
$\omega$ fixed by the group $G_{F|k({\mathbb P}^M)}$ for an appropriate 
$M\ge 1$ and an embedding $k({\mathbb P}^M)\hookrightarrow F$. 
The finite field extension $F^{\ker\rho}|F^K$ can be considered as a purely 
transcendental extension of a function field extension $k(Y)|k(Y)^Q$ of 
smooth projective $k$-varieties of dimension $\ge M$. Consider $\omega$ as 
a differential form with poles on ${\mathbb P}^M_k$. Fix a sufficiently 
general finite morphism $f:Y\longrightarrow{\mathbb P}^M_k$, unramified 
above the poles of $\omega$, and such that the poles of $f^{\ast}\omega$ 
pass through a fixed point of $Y$, but not through another point of its 
$Q$-orbit. Then, as $Q$ acts freely on the set of `sufficiently general' 
divisors on $Y$, the form $p_{\rho}f^{\ast}\omega$ is non-zero, and thus, 
$p_{\rho}f^{\ast}\omega$ spans a $K$-submodule in $W$ isomorphic to $\rho$. 
\qed 

\vspace{4mm} 

{\it Remark.} \label{vanish-I-S-semilin} The vanishing of ${\mathcal I}_S$ 
on any smooth semilinear representation $V$ of $G$ is evident: Let $L$ 
be the function field of an affine $k$-space embedded into $F$. For any 
$v\in V^{G_{F|L}}$ and any $x\in F$ transcendental over $L$ the vector $xv$ 
belongs to $V^{G_{F|L(x)}}$, so its image ${\mathcal I}_SV$ should be fixed 
by $G$. In particular, the image of $xv$ in ${\mathcal I}_SV$ coincides 
with the image of $2xv$, and thus, $xv$ becomes $0$ in ${\mathcal I}_SV$. 
Such vectors $xv$ span $V$, so ${\mathcal I}_SV=0$.

\begin{cor} \label{log-diff-form} For any $k$-variety $U$ and any 
rational closed form $\eta$ on $U\times{\mathbb A}^1$ there exist an 
affine variety $Y$, dominant morphisms $\pi:Y\to U\times{\mathbb A}^1$, 
$\pi_1,\dots,\pi_m:Y\to{\mathbb A}^N_k$ and rational closed forms 
$\eta_1,\dots,\eta_m$ on ${\mathbb A}^N_k$ and $\eta_0$ on $U$ 
such that $\pi^{\ast}\eta=({\rm pr}_U\circ\pi)^{\ast}\eta_0
+\pi_1^{\ast}\eta_1+\dots+\pi_m^{\ast}\eta_m$. \end{cor}
{\it Proof.} We consider $\eta$ as a section of the sheaf 
$\SOm^{\bullet}_{|k,\text{closed}}:X\mapsto
\Omega^{\bullet}_{k(X)|k,\text{closed}}$ over $U\times{\mathbb A}^1$. 
Proposition \ref{A1-quotient-of-closed} describes the kernel of 
$\SOm^{\bullet}_{|k,\text{closed}}\stackrel{\alpha}{\longrightarrow}
{\mathcal I}\SOm^{\bullet}_{|k,\text{closed}}$ as the ideal generated by the 
exact and the logarithmic differentials. By Proposition \ref{I-S-polulin}, 
applied to $W=F^{\times}/k^{\times}$, ${\mathcal I}_S$ annihilates 
the kernel of $\alpha$. Thus, modulo closed forms coming from 
projective spaces, $\eta$ comes from $U$. \qed 

\vspace{4mm} 

Let $X$ be a smooth proper $k$-variety and $W:={\mathbb Q}[\{k(X)\overk F\}]$ 
be the module of generic 0-cycles on $X$. The space $W^K$ is the image of 
the projector defined by the Haar measure of $K$. As the generators of $W$ 
are generic points of $X$, the space $W^K$ is spanned by the 0-cycles of 
type $p_{\ast}\pi^{\ast}q$ for all diagrams of dominant $k$-morphisms 
$X\stackrel{p}{\leftarrow\hspace{-3mm}\longleftarrow}Y
\stackrel{\pi}{\longrightarrow\hspace{-3mm}\to}{\mathbb P}^N_k$, where 
$\pi$ is generically finite, and all generic points $q\in{\mathbb P}^N(F^K)$. 
(Indeed, for any generic $F$-point $\sigma:k(X)\overk F$ of $X$ the orbit 
$K\sigma$ is finite, so the compositum $L_1$ of the images of the elements 
of $K\sigma$ is finitely generated over $k$. Let $L_0\subset F^K$ be a 
finitely generated and purely transcendental extension of $k$ containing 
$L_1^K$. Let $Y$ be a $K$-equivariant smooth $k$-model of $L_0L_1$. Then 
$p$ and $\pi$ are induced by the inclusions $k(X)\subset k(Y)\supset L_0$.) 
Thus, the module ${\mathcal I}_SW$ is the quotient of $W$ by the 
${\mathbb Q}$-span of 0-cycles of type 
$p_{\ast}\pi^{\ast}q_1-p_{\ast}\pi^{\ast}q_2$ for all dominant $k$-morphisms 
$p:Y\to X$, generically finite $k$-morphisms $\pi:Y\to{\mathbb P}^N_k$ and 
all generic points $q_1,q_2\in{\mathbb P}^N(F)$. 

\begin{lemma} \label{cyclic-1} Let $X$ be a smooth proper curve over $k$ 
of genus $g$. Then the $G$-module $Z^{{\rm rat}}_0(k(X)\otimes_kF):=\ker[
Z_0(k(X)\otimes_kF)\longrightarrow CH_0(X\times_kF)]$ is 
generated by $w_N=\sum^N_{j=1}\sigma_j-\sum^N_{j=1}\tau_j$ for all $N>g$, 
where $(\sigma_1,\dots,\sigma_N;\tau_1,\dots,\tau_N)$ is a generic 
$F$-point of the fibre over $0$ of the map 
$X^N\times_kX^N\stackrel{p_N}{\longrightarrow}\Pic^0X$ 
sending $(x_1,\dots,x_N;y_1,\dots,y_N)$ to the class of 
$\sum^N_{j=1}x_j-\sum^N_{j=1}y_j$. \end{lemma} 
{\it Proof.} Let $\gamma_1,\dots,\gamma_s:k(X)\overk F$ and 
$\delta_1,\dots,\delta_s:k(X)\overk F$ be 
generic points of $X$ such that $\sum^s_{j=1}\gamma_j-\sum^s_{j=1}\delta_j$ 
is the divisor of a rational function on $X_F$. 

We need to show that $\sum^s_{j=1}\gamma_j-\sum^s_{j=1}\delta_j$ 
belongs to the $G$-submodule in $Z^{{\rm rat}}_0(k(X)\otimes_kF)$ 
generated by $w_N$'s. 

There is a collection $\alpha_1,\dots,\alpha_g:k(X)
\overk F$ of generic points of $X$ such 
that the class of $\sum^s_{j=1}\gamma_j+\sum^g_{j=1}\alpha_j$ in 
$\Pic^{s+g}X$ is a generic point. Then there is a collection 
$\xi_1,\dots,\xi_{s+g}:k(X)\overk F$ 
of generic points of $X$ in general position such that 
$\sum^s_{j=1}\gamma_j+\sum^g_{j=1}\alpha_j-\sum^{s+g}_{j=1}\xi_j$ 
is divisor of a rational function on $X_F$ (so the same holds 
also for $\sum^s_{j=1}\delta_j+\sum^g_{j=1}\alpha_j-
\sum^{s+g}_{j=1}\xi_j$). We may, thus, assume that 
$\delta_1,\dots,\delta_s$ are in general position. 

Fix a collection $\{\varkappa_{ij}\}_{1\le i\le g,1\le j\le s}$ 
of generic points of $X$ in general position, also with respect 
to $\gamma_1,\dots,\gamma_s$ and to $\delta_1,\dots,\delta_s$, 
such that the classes of $\gamma_1+\sum^g_{i=1}\varkappa_{i1},
\dots,\gamma_s+\sum^g_{i=1}\varkappa_{is}$ in $\Pic^{g+1}X$ 
are generic points in general position. Then one can choose 
a collection $\{\xi_{ij}\}_{0\le i\le g,1\le j\le s}$ of generic 
points of $X$ in general position such that 
$\gamma_j+\sum^g_{i=1}\varkappa_{ij}-\sum^g_{i=0}\xi_{ij}$ 
is divisor of a rational function on $X_F$ (so the same holds 
also for $\sum^s_{j=1}\sum^g_{i=0}\xi_{ij}-\left(\sum^s_{j=1}
\delta_j+\sum^s_{j=1}\sum^g_{i=1}\varkappa_{ij}\right)$). 
We may, thus, assume that both $\gamma_1,\dots,\gamma_s$ 
and $\delta_1,\dots,\delta_s$ are in general position. 

Then there is a collection of generic points $\xi_1,\dots,\xi_s:
k(X)\overk F$ such that the points 
$(\gamma_1,\dots,\gamma_s;\xi_1,\dots,\xi_s)$ and 
$(\delta_1,\dots,\delta_s;\xi_1,\dots,\xi_s)$ are generic 
on $p^{-1}_s(0)$. Then $\sum^s_{j=1}\gamma_j-\sum^s_{j=1}\xi_j$ 
and $\sum^s_{j=1}\delta_j-\sum^s_{j=1}\xi_j$ are divisors 
of rational functions on $X_F$. Clearly, such elements 
belong to the $G$-orbit of $w_s$. \qed 

\vspace{4mm}

{\it Remark.} The $G$-module $Z^{{\rm rat}}_0(k(X)\otimes_kF)$ from 
Lemma \ref{cyclic-1} is generated by $w_{g+1}$. {\it Proof.} There exists 
an effective divisor $D$ (of degree $g$) in the linear equivalence class of 
$\sum^N_{j=2}\sigma_j-\sum^N_{j=g+2}\tau_j$, so $w_N=
[\sum^N_{j=2}\sigma_j-D-\sum^N_{j=g+2}\tau_j]+
[\sigma_1+D-\sum^{g+1}_{j=1}\tau_j]$ is a sum of an element in the 
$G$-orbit of $w_{N-1}$ and an element in the $G$-orbit of $w_{g+1}$. \qed 

\begin{proposition} \label{pic-as-adj} 
${\mathcal I}_S{\mathbb Q}[\{k(X)\overk F\}]=
\Pic(X_F)_{{\mathbb Q}}$ for any smooth proper curve $X$ over $k$. 
\end{proposition}
{\it Proof.} By Lemma \ref{cyclic-1}, it suffices to show that the images 
of the generators $w_N$ in ${\mathcal I}_SW$ are zero. Denote by $g\ge 0$ 
the genus of $X$, by $\psi_N$ a generic effective divisor on $X$ of degree 
$N$ with a special class in $\Pic^NX$. Then $w_N=\sigma\psi_{N+g}-
\tau\psi_{N+g}$ for some $\sigma,\tau\in G$, so it suffices to show that 
the images of $\psi_N$'s in ${\mathcal I}_SW$ are fixed by $G$. 
Denote by $X^N\stackrel{s}{\longrightarrow}\Sigma^NX\stackrel{r}
{\longrightarrow}\Pic^N(X)$ the natural morphisms and set 
$Y:=(rs)^{-1}(\ast)$. Let $p:Y\subseteq X^N\longrightarrow X$ be the 
projection to the first multiple; set $\pi=s|_Y:Y\longrightarrow r^{-1}(\ast)$. 
The projection to the first $N-g$ multiples $Y\longrightarrow X^{N-g}$ 
is generically finite of degree $g!$. If $N\ge 2g-1$ then 
$r^{-1}(\ast)\cong{\mathbb P}^{N-g}$. Assume also that $N\ge g+1$ 
(i.e. $N\ge\max(2g-1,g+1)$). As $s$ is generically finite 
of degree $N!$, one has $(N-1)!\psi_N=p_{\ast}\pi^{\ast}q$ 
for a generic point $q$ of $r^{-1}(\ast)$. \qed 

\vspace{4mm}

Denote by ${\mathcal I}_G$ the full subcategory in $\Sm_G$ of 
``homotopy invariant'' representations: $W^{G_{F|L'}}=W^{G_{F|L}}$ for 
any purely transcendental subextension $L'|L$ in $F|k$. 
\begin{theorem} \label{irr-subq-hom-inv} 
A dominant sheaf is ${\mathbb A}^1$-invariant if and only if all its 
simple subquotients are. \end{theorem}
{\it Proof.} Let $S$ be the collection of pairs of type $(G_{F|L(x)}\subset 
G_{F|L})$ for all subfields $L$ in $F|k$ of finite type and elements 
$x\in F$ transcendental over $L$. The following conditions are equivalent: 
\begin{enumerate} \item \label{poopr} a smooth representation $W$ of $G$ is 
``homotopy invariant''; \item \label{S-dlq-kon-por} $W^{U_1}=W^{U_2}$ 
for all pairs $(U_1\subset U_2)\in S$; \item \label{S-komp} 
$W^{G_{F|L}}=W^{G_{F|L'}}$ for all subfields $L$ in $F|k$ of finite type 
and purely transcendental extensions $L'|L$ in $F$ such that $F$ is 
algebraic over $L'$. \end{enumerate} 
(\ref{poopr})$\Leftrightarrow$(\ref{S-komp}) and 
(\ref{poopr})$\Leftrightarrow$(\ref{S-dlq-kon-por}) are evident; 
(\ref{S-dlq-kon-por})$\Leftrightarrow$(\ref{poopr}) is proved in 
\cite[Corollary 6.2]{repr}. This verifies the condition (2) of 
Proposition \ref{harakteriz-S-objektov}. For each pair $(G_{F|L(x)}\subset 
G_{F|L})\in S$ fix some $\sigma\in G_{F|L}$ with $x$ and $\sigma x$ 
algebraically independent over $L$. Then the condition (1)(i) is obvious: 
$G_{F|L(x)}\cap G_{F|L(\sigma x)}=G_{F|L(x,\sigma x)}$ and 
$(G_{F|L(x,\sigma x)}\subset G_{F|L(x)})\in S$; the condition (1)(ii) 
follows from \cite[Lemma 2.16]{repr}: the subgroups $G_{F|L(x)}$ and 
$G_{F|L(\sigma x)}$ generate $G_{F|L}$. \qed

\subsection{Summary of equivalences} \label{Summ-equiv} 
The following categories are equivalent: 
\begin{enumerate} \item the category of dominant ${\mathbb A}^1$-invariant 
sheaves of $E$-vector spaces; \item the category of dominant 
${\mathbb A}^1$-invariant presheaves of $E$-vector spaces with 
the Galois descent property; 
\item the category $\Sm_G^S(E)$, where $S$ consists of the pairs of type 
$(G_{F|L'}\subset G_{F|L})$ with purely transcendental $L'|L$ in $F|k$. 
\end{enumerate}
These equivalences restrict to equivalences of corresponding subcategories: 
(1) of sheaves of finite-dimensional spaces, (2) of presheaves of 
finite-dimensional spaces, 
(3) of admissible representations of $G$.\footnote{A representation of a 
totally disconnected group is {\sl admissible} if it is smooth and the fixed 
subspaces of all open subgroups are finite-dimensional.} 

Consider the following properties of a smooth representation $W$ of $G$: 
\begin{enumerate} \item \label{A1-inv-predst} $W\in\Sm_G^S(E)$, where $S$ 
consists of the pairs of type $(G_{F|L'}\subset G_{F|L})$ with purely 
transcendental $L'|L$ in $F|k$; 
\item \label{ne-sod-vse} the restriction of $W$ to a compact subgroup 
$U$ does not contain all smooth irreducible representations of $U$; 
\item \label{nenul-ann} the annihilator of $W$ in the algebra 
${\mathbb D}_{{\mathbb Q}}(G)$ is non-zero. \end{enumerate} 
One has (\ref{A1-inv-predst})$\Rightarrow$(\ref{ne-sod-vse})$\Rightarrow$%
(\ref{nenul-ann}). 
[(\ref{ne-sod-vse})$\Rightarrow$(\ref{nenul-ann}) is explained 
in \S\ref{nul-annul-sod-vse}. 
(\ref{A1-inv-predst})$\Rightarrow$(\ref{ne-sod-vse}): 
If $F^U$ is purely transcendental over $k$, there are many irreducible 
smooth representations of $U$, entering in no object of $\Sm_G^S(E)$. 
Any non-trivial smooth irreducible representation $\tau$ of $U$ such that 
$F^{\ker\tau}$ is unirational (e.g., purely transcendental) over $k$ is an 
example of such representation. Clearly, for any such $\tau$ the natural 
projector $p_{\tau}\in{\mathbb D}_{{\mathbb Q}}(G)$ onto the $\tau$-isotypical 
part belongs to the common annihilator of the objects of $\Sm_G^S(E)$.] 

{\it Remark.} \label{glok-por} For a discrete valuation $v$ of rank $1$ on 
$F$, trivial on $k$,\footnote{By definition, this means that any maximal 
system of elements of $F^{\times}$ with independent images in the valuation 
group, should be a transcendence base of $F$ over a lift of a subfield of 
the residue field.} and a smooth representation $W$ of $G$ set 
$W_v:=\sum_LW^{G_{F|L}}\subseteq W$, where $L$ runs over the subfields in 
the valuation ring of $v$. The intersection $\Gamma(W):=\bigcap_vW_v$ over 
all such $v$'s is again in $\Sm_G$. As shown in \cite[Cor.4.7]{max}, the 
property (\ref{A1-inv-predst}) for $W$ implies that $W=W_v$ (and also 
$W=\Gamma(W)$, since all $v$'s as above form a $G$-orbit, cf. \cite{max}). 

\section{Differential forms} \label{Differential-forms}
Let $H^{\bullet}=\bigoplus_{q\ge 0}H^q$ be a cohomology theory, considered as 
a dominant ${\mathbb A}^1$-presheaf. Denote by $\underline{H}_c^{\bullet}$ 
the dominant ${\mathbb A}^1$-sheaf $X\mapsto H^{\bullet}(X)/N^1$ for smooth 
proper $k$-varieties $X$, which is a subsheaf of $\underline{H}^{\bullet}$, 
e.g., $\underline{H}^1_c:X\mapsto H^1(\overline{X})$. Clearly, 
$\underline{H}^{\bullet}_c$ is a sheaf of finite $H^{\bullet}(k)$-modules. 
It would follow from the standard semisimplicity conjecture that the sheaf 
$\underline{H}^{\bullet}_c$ is semisimple if $H^{\bullet}(k)$ is a field. 

We shall be interested in the case of de Rham cohomology $H^{\bullet}
=H^{\bullet}_{{\rm dR}/k}:X\mapsto H^{\bullet}_{{\rm dR}/k}(X):=
{\mathbb H}^{\bullet}(X,\Omega^{\bullet}_{X|k})$, where $H^{\bullet}(k)=k$, 
cf. \cite{Grothendieck}. Clearly, $\underline{H}^q_{{\rm dR}/k}=
\SOm^q_{|k,\text{{\rm closed}}}/\SOm^q_{|k,\text{{\rm exact}}}$, where 
$\Omega^q_{|k,\text{{\rm closed}}}:Y\mapsto
\ker(d|\Gamma(Y,\Omega^q_{Y|k}))$ and $\Omega^q_{|k,\text{{\rm exact}}}:
Y\mapsto d\Gamma(Y,\Omega^{q-1}_{Y|k})$, so $d:
{\mathcal H}^{{\mathbb G}_a}_1\stackrel{\sim}{\longrightarrow}
\SOm^1_{|k,\text{{\rm exact}}}$. The sheaf 
$\underline{H}^1_{{\rm dR}/k,c}$ is semisimple. 
It is described in Lemma \ref{example-4}. 

\subsection{Maximal ${\mathbb A}^1$-subsheaf and the ${\mathbb A}^1$-quotient 
of (closed) forms} 
Recall (\S\ref{Properties-of}) that the inclusion functor 
${\mathcal I}_G\to\Sm_G$ admits a right adjoint $W\mapsto W^{(0)}$, 
the maximal subobject in ${\mathcal I}_G$. The following fact points out 
once more the cohomological nature of the objects of ${\mathcal I}_G$. 
\begin{proposition}[\cite{pgl}, Prop.7.6] \label{wedge} The maximal 
subobject in ${\mathcal I}_G$ of the sheafification of 
$\bigotimes^{\bullet}_{{\mathcal O}}\Omega^1_{|k}$ 
is $\Omega^{\bullet}_{|k,\text{{\rm reg}}}$. 

For any smooth proper $k$-variety $Y$ there are the following 
canonical isomorphisms \begin{equation} \label{ch-0-forms} 
\Hom_{{\mathcal I}_G}(C_{k(Y)},\Omega^q_{|k,\text{{\rm reg}}})=
\Hom_{\Sm_G}(\Psi_Y,\Omega^q_{|k,\text{{\rm reg}}})
\stackrel{\sim}{\longleftarrow}\Gamma(Y,\Omega^q_{Y|k})
\stackrel{\sim}{\longrightarrow}\Hom_{\Sm_G}(CH_0(Y_F),
\Omega^q_{|k,\text{{\rm reg}}}).\end{equation} 
The first isomorphism is functorial with respect to the dominant morphisms 
$Y\longrightarrow Y'$, the second one is functorial with respect to 
arbitrary morphisms $Y\longrightarrow Y'$. \end{proposition}

\begin{lemma} \label{reduction-exact} Let $L$ be an algebraically closed
extension of $k$ and $x$ be an indeterminant. Then there are isomorphisms
$id+\sum_{\alpha\in L}\wedge\frac{d(x-\alpha)}{x-\alpha}:(L(x)
\otimes_L\Omega^q_{L|k})/\Omega^q_{L|k,\text{{\rm exact}}}\oplus\bigoplus
_{\alpha\in L}(\Omega^{q-1}_{L|k}/\Omega^{q-1}_{L|k,\text{{\rm exact}}})
\stackrel{\sim}{\longrightarrow}\Omega^q_{L(x)|k}/
\Omega^q_{L(x)|k,\text{{\rm exact}}}$ and 
$d+\sum_{\alpha\in L}\wedge\frac{d(x-\alpha)}{x-\alpha}:(L(x)
\otimes_L\Omega^q_{L|k})/\Omega^q_{L|k,\text{{\rm closed}}}\oplus\bigoplus
_{\alpha\in L}\Omega^q_{L|k,\text{{\rm exact}}}\stackrel{\sim}{\longrightarrow}
\Omega^{q+1}_{L(x)|k,\text{{\rm exact}}}$ for any $q\ge 1$. The former 
isomorphism restricts to an isomorphism $H^q_{{\rm dR}/k}(L)\oplus
\bigoplus_{\alpha\in L}H^{q-1}_{{\rm dR}/k}(L)
\stackrel{\sim}{\longrightarrow}H^q_{{\rm dR}/k}(L(x))$. \end{lemma}
{\it Proof.} As $\Omega^q_{L(x)|k}=L(x)\otimes_L\Omega^q_{L|k}\oplus
L(x)\otimes_L\Omega^{q-1}_{L|k}\wedge dx$, for any $\omega\in
\Omega^q_{L(x)|k}$ one has $\omega\equiv\eta\wedge dx
\pmod{L(x)\otimes_L\Omega^q_{L|k}}$ for a unique 
$\eta\in L(x)\otimes_L\Omega^{q-1}_{L|k}$. Using partial fraction 
decomposition of rational functions in $L(x)$, we get a presentation
$\eta=\sum_{j\ge 0}x^j\eta_j+\sum_{\alpha\in L,~j\ge 1}
\frac{\eta_{j,\alpha}}{(x-\alpha)^j}$, where
$\eta_j,\eta_{j,\alpha}\in\Omega^{q-1}_{L|k}$. Then $\eta\wedge dx\equiv
\sum_{\alpha\in L}\eta_{1,\alpha}\wedge\frac{d(x-\alpha)}{x-\alpha}
\pmod{L(x)\otimes_L\Omega^q_{L|k}+\Omega^q_{L(x)|k,\text{{\rm exact}}}}$,
so $\omega\equiv\sum_i\phi_i(x)\eta_i+\sum_{\alpha\in L}\eta_{1,\alpha}\wedge
\frac{d(x-\alpha)}{x-\alpha}\pmod{\Omega^q_{L(x)|k,\text{{\rm exact}}}}$,
and thus, $d\omega=\sum_id\phi_i(x)\wedge\eta_i+\sum_i\phi_i(x)d\eta_i+
\sum_{\alpha\in L}d\eta_{1,\alpha}\wedge\frac{d(x-\alpha)}{x-\alpha}\equiv
\sum_i\phi'_i(x)dx\wedge\eta_i+\sum_{\alpha\in L}d\eta_{1,\alpha}\wedge
\frac{dx}{x-\alpha}\pmod{L(x)\otimes_L\Omega^q_{L|k}}$ for some
$\phi_i(x)\in L(x)$ and $\eta_i\in\Omega^q_{L|k}$ (and we may assume that
$\eta_i$ are $L$-linearly independent). Using partial fraction decomposition
of the rational functions $\phi_i\in L(x)$, we see that if $\omega$ is
closed then $d\eta_{1,\alpha}=0$, $\phi_i\in L$ and
$\sum_i\phi_i\eta_i\in\Omega^q_{L|k}$ is closed. \qed

\begin{proposition} \label{A1-quotient-of-closed} Let $M_q$ be the 
sheaf associated with the presheaf $\Omega^q_{|k,\text{{\rm exact}}}+
\Omega^{q-1}_{|k,\text{{\rm closed}}}\wedge d\log{\mathbb G}_m\subset
\Omega^q_{|k}$ for any $q\ge 1$. Then {\rm (i)} 
$\Omega^q_{k(X\times{\mathbb A}^n)|k,\text{{\rm closed}}}
=\Omega^q_{k(X)|k,\text{{\rm closed}}}+M_q(X\times{\mathbb A}^n)$ for any 
$n\ge 1$; {\rm (ii)} $M_q$ is the kernel of the natural projection $\pi_q:
\SOm^q_{|k,\text{{\rm closed}}}\to V^q:={\mathcal I}(\SOm^q_{|k,
\text{{\rm closed}}})={\mathcal I}(\underline{H}^q_{{\rm dR}|k})$; 
{\rm (iii)} for $q\ge 2$, $M_q$ is the sheaf associated with the presheaf 
$\Omega^{q-1}_{|k,\text{{\rm closed}}}\wedge d\log{\mathbb G}_m$ 
and $d+d\log:{\mathcal H}^{{\mathbb G}_a}_1\oplus 
k\otimes{\mathcal H}^{{\mathbb G}_m}_1\to M_1$ is an isomorphism. 
\end{proposition} 
In particular, the natural projections $\SOm^{\bullet}_{|k,\text{{\rm closed}}}
\stackrel{p_1}{\to}\underline{H}^{\bullet}_{{\rm dR}|k}\stackrel{p_2}{\to}
V^{\bullet}:={\mathcal I}(\underline{H}^{\bullet}_{{\rm dR}|k})
={\mathcal I}(\SOm^{\bullet}_{|k,\text{{\rm closed}}})$ are morphisms of 
sheaves of supercommutative $k$-algebras. (The kernel of $p_1$, i.e. 
$\SOm^{\bullet}_{|k,\text{{\rm exact}}}$, is the ideal generated by 
$\SOm^1_{|k,\text{{\rm exact}}}$,\footnote{More generally, let 
$\omega\in\Omega^{\ge i}_{F|k}\smallsetminus\Omega^{\ge i+1}_{F|k}$ 
be a closed form for some $i\ge 0$. Then any ideal in 
$\SOm^{\bullet}_{|k,\text{{\rm closed}}}$ containing the $G$-orbit 
of $\omega$ contains $\SOm^{\ge i+1}_{|k,\text{{\rm exact}}}$. {\it Proof.} 
By \cite[Lemma 7.7]{pgl}, the semilinear representation $\Omega^j_{F|k}$ 
is irreducible for any $j\ge 0$. In particular, $F$-linear envelope of the 
$G$-orbit of $\omega$ is the direct sum of $\Omega^j_{F|k}$ over all $j\ge i$ 
such that the homogeneous component of $\omega$ of degree $j$ is non-zero. 
Then $dz\wedge\sigma\omega=d(z\cdot\sigma\omega)$ for all $z\in F$ and all 
$\sigma\in G$ span the direct sum of $\SOm^j_{|k,\text{{\rm exact}}}$ 
over all $j\ge i$ as above. \qed} the kernel of $p_2$ 
is the ideal generated by $d\log{\mathbb G}_m$.) 
They are surjective even as morphisms of presheaves. 

{\it Proof.} 
Let us show that $\ker\pi_q$ contains $M_q$. For any irreducible smooth 
$k$-variety $X$, any $\eta\in\Omega^{q-1}_{k(X)|k,\text{{\rm closed}}}$ 
and a generator $t$ of the field $k(X\times{\mathbb G}_m)$ over $k(X)$ the 
closed $q$-forms $\omega_m=\eta\wedge d\log t$ and $\omega_a=\eta\wedge dt$ 
are sections of the sheaf $\Omega^q_{|k,\text{{\rm closed}}}$ over 
$X\times{\mathbb G}_m$, so their images in 
${\mathcal I}(\Omega^q_{|k,\text{{\rm closed}}})$ 
should be sections over $X$. As there are endomorphisms $g_m,g_a$ of 
$X\times{\mathbb G}_m|X$ such that $g_mt=t^2$ and $g_at=2t$ (so 
$g_?\omega_?=2\omega_?$), the images of $\omega_?$ in 
${\mathcal I}(\Omega^q_{|k,\text{{\rm closed}}})$ should be zero. The elements 
of type $\eta\otimes d\log t$ (resp., $\eta\otimes dt$) span the sheaf 
$\Omega^{q-1}_{|k,\text{{\rm closed}}}\otimes{\mathcal H}^{{\mathbb G}_m}_1$ 
(resp., $\Omega^{q-1}_{|k,\text{{\rm closed}}}\otimes
{\mathcal H}^{{\mathbb G}_a}_1$, which is surjective over 
$\Omega^q_{|k,\text{{\rm exact}}}$). 

By \cite[Lemma 6.3, p.200]{repr}, to show that $\ker\pi_q=M_q$ 
it suffices to check that, for any algebraically 
closed extension $F'|k$ in $F$ and any $t\in F\smallsetminus F'$, 
any $\omega\in\Omega^q_{F'(t)|k,\text{{\rm closed}}}$ 
belongs in fact to $\Omega^q_{F'|k,\text{{\rm closed}}}+M_q$. 

By Lemma \ref{reduction-exact}, $\omega\equiv 
\xi+\sum_{\alpha\in F'}\eta_{\alpha}\wedge\frac{d(t-\alpha)}{t-\alpha}
\pmod{\Omega^q_{F'(t)|k,\text{{\rm exact}}}}$, where 
$\xi\in\Omega^q_{F'|k,\text{{\rm closed}}}$ and 
$\eta_{\alpha}\in\Omega^{q-1}_{F'|k,\text{{\rm closed}}}$, which 
means that $\omega\in\Omega^q_{F'|k,\text{{\rm closed}}}+M_q$. \qed 

\begin{conjecture} The sheaf $V^{\bullet}$ is semisimple. \end{conjecture}

\vspace{4mm}

{\it Remarks.} 1. It follows from Proposition \ref{A1-quotient-of-closed} 
that the natural morphism $\underline{H}^{\bullet}_{{\rm dR}/k,c}\to 
V^{\bullet}$ is injective. 

2. As explained in Remark on p.\pageref{vanish-I-S-semilin}, 
${\mathcal I}V=0$ for any semilinear smooth representation $V$: 
if $v\in V^{G_{F|L}}$ and $f\in F$ is transcendental over $L$ then 
$v=fv-(f-1)v$ becomes zero in any quotient of $V$ in ${\mathcal I}_G$.

3. For an algebra $A\in\Sm_G$ it is not always true that the kernel 
$A^{\circ}$ of the projection $A\to{\mathcal I}A$ is an ideal. E.g., let 
$A=A_{\bullet}$ be the (graded) tensor, symmetric or skew-symmetric algebra 
of $A_1={\mathbb Q}[F\smallsetminus k]$. Then ${\mathcal I}A_1={\mathbb Q}$, 
so $A_1^{\circ}\otimes A_1+A_1\otimes A_1^{\circ}$ consists of all 
sums in ${\mathbb Q}[(F\smallsetminus k)\times(F\smallsetminus k)]$ of degree 
0. On the other hand, ${\mathcal I}(A_1\otimes A_1)=\bigoplus_{x\in
\mathop{\mathbf{Spec}}(k({\mathbb P}^1)\otimes_kk({\mathbb P}^1))}C_{k(x)}$, 
and therefore, $A_1^{\circ}\otimes A_1+A_1\otimes A_1^{\circ}$ is 
strictly bigger than $A_2^{\circ}$. 

\subsection{The semisimplicity of the regular forms of top degree} 
\label{semisimpl-top-sect} 
Let $L$ be an algebraically closed extension of $k$ with 
$1\le q=\mathop{\mathrm{tr.deg}}(L|k)<\infty$. Define a representation 
$\Omega^q_{L|k,\text{{\rm reg}}}$ as the union in $\Omega^q_{L|k}$ of 
all spaces $\Gamma(X,\Omega^q_{X|k})$ over all smooth proper varieties 
$X$ over $k$ with the function field embedded into $L$ over $k$. 

The na\"{\i}ve truncation filtration on $\Omega^{\bullet}_{\overline{U}|k}$ 
gives the descending Hodge filtration $F^{\bullet}$ on 
$H^q_{{\rm dR}/k}(\overline{U})$. The Hodge filtrations on 
$H^q_{{\rm dR}/k}(\overline{U})$ for all $U$'s induce a canonical 
filtration $F^{\bullet}$ on $\underline{H}^q_{{\rm dR}/k,c}$ 
by subsheaves of $k$-vector spaces with associated graded quotients 
$H^{p,q-p}_{|k}:Y\mapsto{\rm coker}[\bigoplus_D
H^{p-1}(D,\Omega^{q-p-1}_{D|k})\longrightarrow 
H^p(\overline{Y},\Omega^{q-p}_{\overline{Y}|k})]$, where 
$D\to\overline{Y}$ runs over all resolutions of the divisors on 
$\overline{Y}$. In particular,  
$H^{q,0}_{|k}=F^q\underline{H}^q_{{\rm dR}/k,c}=\Omega^q_{|k,\text{reg}}:
Y\mapsto\Gamma(\overline{Y},\Omega^q_{\overline{Y}|k})$ is the dominant 
subsheaf of $\underline{H}^q_{{\rm dR}/k,c}$ consisting of regular 
differential $q$-forms. 

\begin{proposition}[\cite{glatt}] \label{semisimpl-top} 
Suppose that the cardinality of $k$ is at most continuum. 
The representation $H^q_{{\rm dR}/k,c}(L)$ $($and therefore, 
$\Omega^q_{L|k,\text{{\rm reg}}})$ of $G_{L|k}$ is semisimple. Any 
embedding $\iota:k\hookrightarrow{\mathbb C}$ into the field of complex 
numbers determines \begin{itemize} \item a ${\mathbb C}$-antilinear 
isomorphism $H^{s,t}_{L|k}\otimes_{k,\iota}{\mathbb C}\cong 
H^{t,s}_{L|k}\otimes_{k,\iota}{\mathbb C}$, \item a positive definite 
$G_{L|k}$-equivariant hermitian form $({\mathbb C}\otimes_{k,\iota}
H^q_{{\rm dR}/k,c}(L))\otimes_{id,{\mathbb C},\sigma}({\mathbb C}
\otimes_{k,\iota}H^q_{{\rm dR}/k,c}(L))\longrightarrow
{\mathbb C}(\chi)$, where $\sigma$ is the complex conjugation and 
$\chi$ is the modulus of $G_{L|k}$. \end{itemize} There exists a 
non-canonical ${\mathbb Q}$-linear isomorphism $H^{s,t}_{L|k}\cong 
H^{t,s}_{L|k}$. \end{proposition}
{\it Proof.} For any smooth projective $k$-variety $X$ the complexified 
projection $F^pH^{p+q}_{{\rm dR}/k}(X)\to H^q(X,\Omega^p_{X|k})$ 
identifies $F^pH^{p+q}_{{\rm dR}/k}(X)\otimes_{k,\iota}{\mathbb C}\cap
\overline{F^qH^{p+q}_{{\rm dR}/k}(X)\otimes_{k,\iota}{\mathbb C}}$ 
with $H^q(X,\Omega^p_{X|k})\otimes_{k,\iota}{\mathbb C}$. This gives 
a decomposition ${\mathbb C}\otimes_{k,\iota}H^q_{{\rm dR}/k,c}(L)=
\bigoplus_{s+t=q}{\mathbb C}\otimes_{k,\iota}H^{s,t}_{L|k}$. Then the 
complex conjugation on $H^{p+q}(X_{\iota}({\mathbb C}),{\mathbb C})=
H^{p+q}(X_{\iota}({\mathbb C}),{\mathbb R})\otimes_{{\mathbb R}}{\mathbb C}$ 
identifies $H^q(X,\Omega^p_{X|k})\otimes_{k,\iota}{\mathbb C}$ with 
$H^p(X_{\iota}({\mathbb C}),\Omega^q_{X_{\iota}({\mathbb C})})=
H^p(X,\Omega^q_{X|k})\otimes_{k,\iota}{\mathbb C}$. 

The semisimplicity of the $k$-representation $H^q_{{\rm dR}/k,c}(L)$ 
of $G_{L|k}$ is equivalent to the semisimplicity of its complexification. 
For the latter note that there is a positive 
definite $G_{L|k}$-equivariant hermitian form $({\mathbb C}\otimes_{k,\iota}
H^{s,t}_{L|k})\otimes_{id,{\mathbb C},\sigma}({\mathbb C}\otimes_{k,\iota}
H^{s,t}_{L|k})\longrightarrow{\mathbb C}(\chi)$, given by $(\omega,\eta)=
\int_{X_{\iota}({\mathbb C})}i^{q^2+2t}\omega\wedge\overline{\eta}\cdot
[G_{L|k(X)}]$ for any $\omega,\eta\in H^{s,t}_{{\rm prim}}(X_{\iota}
({\mathbb C}))=H^t_{{\rm prim}}(X,\Omega^s_{X|k})\otimes_{k,\iota}
{\mathbb C}\subset{\mathbb C}\otimes_{k,\iota}H^{s,t}_{L|k}$. 
Here $H^{s,t}_{{\rm prim}}(X_{\iota}({\mathbb C}))$ denotes the subspace 
orthogonal to the sum of all Gysin maps $H^{s-1,t-1}(D)\longrightarrow 
H^{s,t}(X_{\iota}({\mathbb C}))$ for all desingularizations $D$ of all 
divisors on $X_{\iota}({\mathbb C})$, as in the definition of 
$\Omega^q_{L|k,\text{{\rm reg}}}$, $X$ runs over all smooth proper 
$k$-varieties with the function field embedded into $L|k$. \qed 

\subsection{Structure of closed 1-forms}
Let ${\rm Div}^{\circ}_{{\mathbb Q}}:Y\mapsto{\rm Div}_{{\rm alg}}
(\overline{Y})_{{\mathbb Q}}$ be the presheaf of algebraically trivial 
divisors. It is a sheaf. 
\begin{lemma} \label{surj-res} The residue homomorphism ${\rm Res}_Y:
H^1_{{\rm dR}/k}(k(Y))\to k\otimes{\rm Div}(\overline{Y})$, 
$\omega\mapsto({\rm res}_x\omega)_{x\in\overline{Y}^1}$, defines 
a morphism of sheaves ${\rm Res}:\underline{H}^1_{{\rm dR}/k}\to 
k\otimes{\rm Div}^{\circ}_{{\mathbb Q}}$. The short sequence 
$0\to\underline{H}^1_{{\rm dR}/k,c}\to\underline{H}^1_{{\rm dR}/k}
\stackrel{{\rm Res}}{\longrightarrow}{\rm Div}^{\circ}_{{\mathbb Q}}
\otimes k\to 0$ is exact, even as a sequence of presheaves. \end{lemma} 
{\it Proof.} As ${\rm Res}$ commutes with the restriction to any 
sufficiently general curve $C$,
${\rm Res}_X(\omega)\cdot C={\rm Res}_C(\omega |_C)\in CH_0(X)$,
$\deg({\rm Res}_X(\omega)\cdot C)=0$ by Cauchy theorem,
the pairing ${\rm NS}(X)_{{\mathbb Q}}\otimes CH_1(X)_{{\mathbb Q}}/hom
\longrightarrow{\mathbb Q}$ is non-degenerate (by Lefschetz hyperplane
section theorem), the class of ${\rm Res}_X(\omega)$ in
${\rm NS}(X)_{{\mathbb Q}}$ is zero. Thus, ${\rm Res}_X$
factors through  the algebraically trivial divisors on $X$.

Clearly, the kernel of ${\rm Res}$ coincides with $\underline{H}^1
_{{\rm dR}/k,c}$, cf. \cite{rosenlicht}.\footnote{If the residues of 
$\omega\in H^1_{{\rm dR}/k}(k(X))$ are zero then integration along a loop 
depends only on its homology class in $H_1(X,{\mathbb Q})$. There is an 
element $\eta$ of $H_{{\rm dR}/k}^1(X)$ with the same periods as $\omega$, 
so integration of $\omega-\eta$ along a path joining a fixed 
(rational) point with the variable one is independent of a 
chosen path, and defines a meromorphic (i.e. rational) function.} Then 
it remains to show that any algebraically trivial divisor on $X$ is 
the residue of a closed 1-form. Any algebraically trivial divisor can 
be written as $D_1-D_2$ for a pair $D_1,D_2$ of algebraically equivalent 
effective divisors on $X$. There is a smooth projective curve $C$, and 
an effective divisor $D$ on $X\times C$, such that ${\rm pr}_X:D\to X$ 
is generically finite and $D_P-D_Q=D_1-D_2$ for some points $P,Q\in C$. 
By Riemann--Roch theorem for curves, there exists a 1-form 
$\omega_{P,Q}\in\Omega^1_C(P+Q)$ such that ${\rm Res}_C(\omega_{P,Q})=P-Q$: 
there is a non-holomorphic 1-form with simple poles in the set $\{P,Q\}$, 
since $\dim_k\Gamma(C,\Omega^1_C(P+Q))=\dim_k\Gamma(C,\Omega^1_C)+1$; 
there are no 1-forms with precisely one simple pole, since 
$\Gamma(C,\Omega^1_C(P))=\Gamma(C,\Omega^1_C(Q))=\Gamma(C,\Omega^1_C)$. 
Then ${\rm Res}_X({\rm pr}_{X\ast}(({\rm pr}_C^{\ast}\omega_{P,Q})|_D))
=D_1-D_2$. \qed 

\begin{proposition} \label{1-forms} 
\begin{itemize} \item The maximal semisimple subsheaf of 
$\SOm^1_{|k,\text{{\rm closed}}}$ is canonically isomorphic to the direct 
sum $\bigoplus_A\Gamma(A,\Omega^1_{A|k})^{A(k)}\otimes_{\End(A)}
{\mathcal H}^A_1={\mathcal H}^{{\mathbb G}_a}_1\oplus 
k\otimes{\mathcal H}^{{\mathbb G}_m}_1\oplus\Omega^1_{|k,\text{{\rm reg}}}$, 
where $A$ runs over the set of isogeny classes of simple commutative 
algebraic $k$-groups; $\Gamma(A,\Omega^1_{A|k})^{A(k)}=\Hom_k({\rm Lie}(A),k)$ 
denotes the space of translation invariant 1-forms on $A$. The projection 
$\SOm^1_{|k,\text{{\rm closed}}}\to\SOm^1_{|k,\text{{\rm closed}}}/
\Omega^1_{|k,\text{{\rm reg}}}$ is split (but not canonically). 
\item The maximal semisimple subsheaf of $\underline{H}^1_{{\rm dR}/k}$ 
is canonically isomorphic to $\bigoplus_AH^1_{{\rm dR}/k}(A)
\otimes_{\End(A)}{\mathcal H}^A_1=k\otimes{\mathcal H}^{{\mathbb G}_m}_1
\oplus\underline{H}^1_{{\rm dR}/k,c}$, where $A$ runs over the set of 
isogeny classes of simple commutative algebraic $k$-groups (with the zero 
summand corresponding to ${\mathbb G}_a$). The projection 
$\underline{H}^1_{{\rm dR}/k}\to\underline{H}^1_{{\rm dR}/k}/
\underline{H}^1_{{\rm dR}/k,c}$ is split (but not canonically). 
\item The sheaf $V^1:Y\mapsto H^1_{{\rm dR}/k}(k(Y))/
k\otimes(k(Y)^{\times}/k^{\times})$ from Proposition 
\ref{A1-quotient-of-closed} is canonically isomorphic to 
$\bigoplus_AV^1(A)\otimes_{\End(A)}{\mathcal H}^A_1$, where $A$ runs 
over the set of isogeny classes of simple abelian $k$-varieties. 
\end{itemize} For any integer $q\ge 1$, the representation 
$\Omega^1_{L|k,\text{{\rm closed}}}$ of the group $G_{L|k}$ 
admits similar description (cf. \S\ref{semisimpl-top-sect}). 
\end{proposition} 
{\it Proof.} In notation of Lemma \ref{surj-res}, the sheaf 
${\rm Div}^{\circ}_{{\mathbb Q}}$ admits a natural surjective 
morphism onto the Picard sheaf ${\rm Pic}^{\circ}_{{\mathbb Q}}=
{\rm coker}[{\mathcal H}^{{\mathbb G}_m}_1
\stackrel{{\rm div}}{\longrightarrow}{\rm Div}^{\circ}_{{\mathbb Q}}]:
Y\mapsto{\rm Pic}^0(\overline{Y})_{{\mathbb Q}}$ with the irreducible 
kernel ${\mathcal H}^{{\mathbb G}_m}_1$. The Picard sheaf 
${\rm Pic}^{\circ}_{{\mathbb Q}}$ is semisimple and 
it is described in Lemma \ref{example-4}. 

According to Lemma \ref{example-4}, for any simple abelian variety $A$ over 
$k$, any non-zero element $\xi$ of ${\rm Pic}^0(A)(k)_{{\mathbb Q}}$ provides 
an embedding of ${\mathcal H}^A_1$ into ${\rm Pic}^{\circ}_{{\mathbb Q}}$. 
Let us show that the natural extension $0\to{\mathcal H}^{{\mathbb G}_m}_1
\to{\rm Div}^{\circ}_{{\mathbb Q}}\to{\rm Pic}^{\circ}_{{\mathbb Q}}\to 0$ 
does not split, even after restricting to ${\mathcal H}^A_1$ via $\xi$. 

All elements of ${\rm Pic}^{\circ}_{{\mathbb Q}}(A):=
{\rm Pic}^0(A)_{{\mathbb Q}}$ are fixed by translations of $A$ by torsion 
elements in $A(k)$. However, as the torsion subgroup in $A(k)$ is Zariski 
dense, it cannot fix a non-zero element of ${\rm Div}^{\circ}_{{\mathbb Q}}(A)
:={\rm Div}_{\text{alg}}(A)_{{\mathbb Q}}$. 

This implies that ${\mathcal H}^{{\mathbb G}_m}_1$ is the maximal 
semisimple subsheaf of ${\rm Div}^{\circ}_{{\mathbb Q}}$, which proves, 
by Lemma \ref{surj-res}, the second assertion. It follows also that the 
simple subquotients of $V^1$ are isomorphic to ${\mathcal H}^A_1$ for 
simple abelian $k$-varieties $A$. There are no extensions between 
${\mathcal H}^A_1$ and ${\mathcal H}^B_1$ for abelian $k$-varieties $A$ 
and $B$, since ${\mathcal I}_G$ is a Serre subcategory of $\Sm_G$ by 
Proposition \ref{harakteriz-S-objektov} and $Y\mapsto A(k(Y))_{{\mathbb Q}}$ 
is a projective object of ${\mathcal I}_G$ by property 5 of 
\S\ref{Properties-of} and Proposition \ref{pic-as-adj}. This 
means that $V^1$ is semisimple, which proves the third assertion. 

Once we know the simple subquotients of $\SOm^1_{|k,\text{{\rm closed}}}$, 
the first assertion follows from Proposition \ref{wedge} and Lemma 
\ref{example-4}. To see that the projections $\SOm^1_{|k,\text{{\rm closed}}}
\to\SOm^1_{|k,\text{{\rm closed}}}/\Omega^1_{|k,\text{{\rm reg}}}$ and 
$\underline{H}^1_{{\rm dR}/k}\to\underline{H}^1_{{\rm dR}/k}/
\underline{H}^1_{{\rm dR}/k,c}$ are split, it is enough to notice that $V^1$ 
is semisimple, and therefore, the compositions $\Omega^1_{|k,\text{{\rm reg}}}
\hookrightarrow\SOm^1_{|k,\text{{\rm closed}}}\to V^1$ and 
$\underline{H}^1_{{\rm dR}/k,c}\hookrightarrow\underline{H}^1_{{\rm dR}/k}\to 
V^1$ admit splittings. \qed

\vspace{4mm}

\label{coniveau-filtr} 
Given a subfield $L$ in $F$, define the filtration $N_{\bullet}^{(L)}$ on 
the $G_{F|L}$-modules $W$ by $N_j^{(L)}W=\sum_{F'}W^{G_{F|F'}}$, where $F'$ 
runs over the subfields in $F|L$ of transcendence degree $j$. (Clearly, 
$N_0^{(L)}W=W^{G_{F|\overline{L}}}$ and $N_{\bullet}^{(L)}=
N_{\bullet}^{(\overline{L})}\subseteq N_{\bullet}^{(L')}$ if 
$L\subset L'\subset F$.) In particular, define the {\sl level} 
filtration on the $G$-modules by $N_{\bullet}:=N_{\bullet}^{(k)}$. 

It is conjectured in \cite[Conj.6.9]{repr} that the graded pieces 
of $N_{\bullet}$ on the objects of ${\mathcal I}_G$ are semisimple. 

If $U$ is an open subgroup of $G$, contained between $G_{F|\overline{L}}$ 
and the normalizer of $G_{F|\overline{L}}$ in $G$ then $N_{\bullet}^{(L)}$ 
is a filtration by $U$-submodules. The forgetful functor $\Sm_G\to\Sm_U$ 
does not preserve the irreducibility (or the semisimplicity). E.g., 
for any commutative simple algebraic $k$-group $A$, the restriction to $U$ 
of the irreducible $G$-module $A(F)/A(k)$ is a non-split extension of 
the irreducible $U$-module $A(F)/A(\overline{L})$ by the $U$- (in fact, 
$(U/U\cap G_{F|\overline{L}})$- ) module $A(\overline{L})/A(k)$. 

{\it Questions.} Let $W\in\Sm_G$ be irreducible 
and $W=N_qW$. Is it true that the representation $W/N_{q-1}^{(L)}W$ of $U$ 
is irreducible (or zero)? 

Clearly, $N_j\Omega^i_{F|k}=\Omega^i_{F|k}$ for any $j>i$, 
$N_j\Omega^i_{F|k}=0$ for any $j<i$ and 
$N_j\Omega^j_{F|k}\subseteq\Omega^j_{F|k,\text{closed}}$. 
\begin{conjecture} $N_j\Omega^j_{F|k}=N_j\Omega^j_{F|k,\text{{\rm closed}}}
=\Omega^j_{F|k,\text{{\rm closed}}}$. \end{conjecture} 
A ``weak'' version, 
$N_j\Omega^j_{F|k,{\rm reg}}=\Omega^j_{F|k,{\rm reg}}$, 
follows from Grothendieck's diagonal decomposition conjecture. 

The Conjecture obviously holds true for $j=0$. The case $j=1$ follows 
from (i) Proposition \ref{1-forms}, (ii) the fact (\cite[Cor.3.8]{glatt}) 
that $F/k$ and $F^{\times}/k^{\times}$ 
are acyclic, so $\Omega^1_{L|k,\text{closed}}\to\hspace{-3mm}\to H^0(G_{F|L},
H^1_{\text{dR}/k}(F)/k d\log(F^{\times}/k^{\times}))$, 
(iii) $N_1(A(F)/A(k))=A(F)/A(k)$ for any commutative $k$-group $A$. \qed

\vspace{5mm}

\noindent
{\sl Acknowledgements.} {\small The project originates from my stay at 
the Max-Planck-Institut in Bonn, it has reached its present state at 
the Institute for Advanced Study in Princeton, and it has reached the 
final form at the I.H.E.S. in Bures-sur-Yvette. I am grateful to these 
institutions for their hospitality and exceptional working conditions. 

}

\end{document}